\documentclass[11pt,twoside,reqno]{article}

\usepackage{amssymb}

\usepackage{amsmath}

\usepackage{amsfonts}
\usepackage{graphicx}
\usepackage{graphics}
\usepackage{amssymb}

\usepackage{color}
 \usepackage{bbm}
\usepackage{nomencl}

\usepackage{titlesec}
\usepackage[titletoc,toc,title]{appendix}

\titleformat{\section}{\large\bfseries}{\thesection.}{0.5em}{}
\titleformat{\subsection}{\bfseries}{\thesubsection.}{0.5em}{}
\titleformat{\subsubsection}{\bfseries}{\thesubsubsection.}{0.5em}{}

\makeglossary
 \allowdisplaybreaks[2]

\usepackage{latexsym}

\DeclareMathOperator{\SL}{SL} 
\DeclareMathOperator{\SO}{SO} \DeclareMathOperator{\Oo}{O}

\renewcommand{\skew}{\mathop{\rm skew}}
 \DeclareMathOperator{\sym}{sym}
\DeclareMathOperator{\tr}{tr} \DeclareMathOperator{\axl}{axl}
\DeclareMathOperator{\anti}{anti} \DeclareMathOperator{\dev}{dev}
\DeclareMathOperator{\Curl}{Curl\hskip.04truecm}
\DeclareMathOperator{\Grad}{Grad\hskip.04truecm}
\DeclareMathOperator{\Div}{Div} \DeclareMathOperator{\curl}{curl\,}
\DeclareMathOperator{\inc}{inc}
\DeclareMathOperator{\sL}{\mathfrak{sl}}
\DeclareMathOperator{\so}{\mathfrak{so}}

\newcommand{\yieldlimit}{\sigma_{\mathrm{y}}}

\newcommand{\yieldzero}{{\sigma}_{0}}

\newcommand{\C}{\mathbb{C}}

\newcommand{\BBR}{\mbox{$\mathbb{R}$}}

\newcommand{\SFH}{\mbox{$\mathsf{H}$}}
\newcommand{\SFL}{\mbox{$\mathsf{L}$}}

\newcommand{\SFW}{\mbox{$\mathsf{W}$}}

\newcommand{\SFZ}{\mbox{$\mathsf{Z}$}}

\newcommand{\bvarepsilon}{\mbox{$\bf\varepsilon$}}

\newcommand{\dsize}{\displaystyle}

\newcommand{\qed}{\qquad$\blacksquare$}
\renewcommand{\div}{\mathop{\rm div}\nolimits}
\newcommand{\ba}{\mbox{\boldmath{$a$}}}
\numberwithin{equation}{section}

\newcommand{\norm}[1]{\left\lVert {#1} \right\rVert} 
\newcommand{\bfig}[2]{\begin{figure}\begin{center}\begin{picture}(341.8,#2)(
#1,0)}
\newcommand{\efig}[2]{\end{picture}\caption{#2.}\lbl{#1}\end{center}
\end{figure}}

\newcommand{\la}{\langle}
\newcommand{\ra}{\rangle}
\newcommand{\id}{{\boldsymbol{\mathbbm{1}}}}

\newtheorem{theorem}{Theorem}[section]

\newtheorem{remark}{Remark}[section]

\newcommand{\be}{\begin{equation}}

\newcommand{\ee}{\end{equation}}


\makeatletter
 \let\@fnsymbol\@arabic
 \makeatother
\begin{document}
\vskip-3truecm 
\title{
\vspace{-1.in} {\Large A fourth order gauge-invariant gradient plasticity model for polycrystals 
 based on Kr\"{o}ner's incompatibility tensor\texttt{}}}

\author{{\large Fran\c{c}ois Ebobisse{\footnote{Corresponding author, Fran\c{c}ois Ebobisse,
Department of Mathematics and Applied Mathematics, University of
Cape Town, Rondebosch 7700, South Africa, e-mail:
francois.ebobissebille@uct.ac.za }}
\quad and\quad Patrizio
Neff{\footnote
{Patrizio Neff, Head of Chair for  Nonlinear Analysis, Faculty of Mathematics, University of
Duisburg-Essen, Thea-Leymann Str. 9, 45127 Essen, Germany, e-mail:
patrizio.neff@uni-due.de, http://www.uni-due.de/mathematik/ag\b{
}neff}}}\vspace{1mm}}



\date{\today}

 \maketitle

\begin{center}

\vspace{-1.5cm}\noindent

{\em Dedicated to Elias C. Aifantis, pioneer of gradient
plasticity}

\vspace{9mm}
 {{\small
\begin{abstract} In this paper we derive a novel fourth order gauge-invariant phenomenological model of infinitesimal
rate-independent gradient plasticity with isotropic hardening  and
Kr\"{o}ner's incompati\-bility tensor
$\inc(\bvarepsilon_p):=\Curl[(\Curl\bvarepsilon_p)^T]$, where $\bvarepsilon_p$ is the symmetric  plastic strain tensor. Here, gauge-invariance denotes invariance under diffeomorphic  reparametrizations of the reference configuration, suitably adapted to the geometrically linear setting.
 The model
features a defect energy contribution which is quadratic in the
tensor $\inc(\bvarepsilon_p)$ and it contains isotropic hardening
based on the rate of the plastic strain
tensor $\dot{\bvarepsilon}_p$.  We motivate the new model by introducing a novel
rotational invariance requirement in gradient plasticity, which we
call micro-randomness, suitable for the description of polycrystalline aggregates on a
mesoscopic scale and not coinciding with classical isotropy
requirements. This new condition effectively reduces the increments of the non-symmetric plastic distortion $\dot{p}$ to their symmetric counterpart $\dot{\bvarepsilon}_p=\sym\dot{p}$.  In the polycrystalline case, this condition is a statement about insensitivity to arbitrary superposed grain rotations. We formulate a mathematical existence result for a suitably regularized  non-gauge-invariant model. The regularized model is rather invariant under reparametrizations of the reference configuration including infinitesimal conformal mappings.
   \end{abstract}}}
\end{center}
\newpage
 {\small \tableofcontents}
\section{Introduction}\label{Intro}
In recent years there has been a growing attention in extending
continuum plasticity theories towards the incorporation of the
experimentally observed size-effects in small scales (see e.g. \cite{Fleckhutch1997, FMAH1994, NIXGAO1998,Stolkevans1998}). This extension
is mainly done via the introduction of certain gradient terms,
making the plastic evolution in some sense nonlocal. Perhaps the
earliest such phenomenological model is due to Aifantis et al.
\cite{AIF1984,AIF1987,MUHAIF1991,ZBIBAIF1992}, who directly
incorporated the Laplacian $\Delta\gamma_p$ in the flow stress,
where $\gamma_p:=\int_0^t\norm{\dot{\bvarepsilon}_p}\,ds$ is a measure of
accumulated equivalent plastic strain\footnote{Cf. e.g. Hill
\cite[p.30]{Hill1950}.} (see Table \ref{table:summary-aif} below for
a summary). \begin{table}[h!]\footnotesize \begin{center}
\begin{tabular}{|l|}\hline \qquad\\

$\sym\nabla u =\varepsilon=\varepsilon_e +\varepsilon_p$\,,\quad
$\varepsilon_p\in\mbox{Sym}\,(3)$\\[.72em]
 $\yieldlimit=\sigma_0+\mu\,k_2\,\gamma_p-\mu\,L_c^2\,\Delta\gamma_p$\,,  \,\,\,$\yieldlimit$ is
the current yield stress\\[.72em]
$L_c>0$ material length scale\\[.72em]
 $\sigma_0>0$\,\, the initial yield stress\\[.72em]
$k_2$\,\, can be taken negative for local softening behaviour
\\[.72em]
$\gamma_p:=\int_0^t\norm{\dot{\varepsilon}_p}\,ds$ \,\,the accumulated equivalent plastic strain\,.
\\\\
\hline  \end{tabular}
 \end{center}\caption{\footnotesize A summary of the model by Aifantis \cite{AIF1984, AIF1987,MUHAIF1991} based on the accumulated equivalent plastic strain $\gamma_p$.
 Depending on the sign of $\Delta\gamma_p$, the model describes process-dependent hardening ($\Delta\gamma_p<0$)
  or softening ($\Delta\gamma_p>0$) due to nonlocal effects. }
 \label{table:summary-aif}\end{table}
 Other variants of the Aifantis' model also based on the accumulated plastic strain were proposed later through
 the principle of virtual power by Fleck and Hutchinson
 \cite{Fleckhutch2001}, Gudmundson \cite{GUD2004} and generalized by Gurtin and Anand \cite{GURTAN2009}.

While there are numerous proposals of such gradient enhanced
phenomenological models, either based on the multiplicative
decomposition
\begin{equation}\label{multi-split1}F=F_e\cdot F_p\end{equation} (see e.g.
\cite{ACHABASS2000,CERGURT2001, GURTAN2005FD, NCA}), or based on the
geometrically linearized corresponding additive decomposition
\begin{equation}\label{add-split1}\bvarepsilon=\bvarepsilon_e+\bvarepsilon_p\end{equation} (see
e.g. \cite{AIF1984, AIF1987,MUHAIF1991,Fleckhutch2001,GUD2004,
GURTAN2005,GURTAN2009,GURED2014,MENZSTEIN1998, MENZSTEIN,MUHAIF1991,SVENEFFMENZ2009}), no general consensus has been reached as to
which variables describing the plastic evolution should be employed
and how they should be combined with their partial derivatives in
space. For example,  a dependence of the plastic flow rule directly
via the infinitesimal plastic strain variable $\bvarepsilon_p$ is sometimes
excluded, since the backstress variable $\bvarepsilon_p$ is not {\bf gauge-invariant}.{\footnote{Aifantis writes in \cite[p. 218]{AIF1987}:  "...
In conformity with established results - that the plastic strain
rate $\dot{\bvarepsilon}_p$ is a state variable, rather than the
strain $\bvarepsilon_p$ itself."}

However, if $\bvarepsilon_p$ is not allowed to appear itself in the
equations, then linear kinematic hardening \'a la Prager is excluded
from the onset, while the modelling of classical linear isotropic
hardening remains possible since it is based on
$\norm{\dot{\bvarepsilon}_p}$, which is invariant under reparametrization of the reference configuration ({\bf gauge-invariance}), see subsection \ref{gauge-inv} in this paper.

With the use of an evolution equation for the symmetric plastic
strain tensor $\bvarepsilon_p$ being traditional, there are also
approaches which focus directly on the plastic distortion $p$ (which
is a non-symmetric variable), thus allowing for the so-called
{\bf plastic spin}{\footnote{The plastic spin in the finite deformation flow 
theory of plasticity is defined as the skew-symmetric part of the
so-called plastic distortion rate i.e.,
$W_p:=\skew(\dot{F}_pF_p^{-1})$\,, while
 its counter-part in the small strain theory is simply $\skew(\dot{p})$\,, where $p$ is the non-symmetric infinitesimal plastic distortion.}} (see
\cite[p.493 and eqs. (91.7) and (91.10)]{GURTAN-BOOK} and
also \cite{Dafalias1985, Kratochvil1973, Mandel1973, Bardella2009,
POH2013, Bardepante2015}). The connection between the two approaches
is simple: we can always identify $\bvarepsilon_p:=\sym p$. As it
will turn out subsequently, the introduction of the plastic
distortion $p$ will make our modelling framework much more
transparent: on the one hand, the passage from the multiplicative
decomposition to the additive decomposition via formal  geometric
linearization is easier and the discussion of invariance conditions
 becomes clearer, even if in the end we obtain a model for the
plastic strain tensor $\bvarepsilon_p=\sym p$ only.

Our aim with this paper is to present a rational modelling environment
for gradient plasticity with respect to the small strain framework
and the additive decomposition which incorporates certain insights
learned from the multiplicative decomposition.

Let us therefore collect what the model should be able to do. It
should
\begin{itemize}
\item[-] incorporate energetic hardening (due to {\bf G}eometrically {\bf N}ecessary
{\bf D}islocations, GNDs) related to the energetic length scale $L_c$;
\item[-] incorporate nonlocal hardening (backstress and Bauschinger effects);
\item[-] satisfy appropriate invariance conditions (objectivity,
referential isotropy, independence of refe\-rence configuration,
gauge-invariance, elastic isotropy, elastic frame-indifference,
etc.);
\item[-] allow in principle for plastic spin (Bardella
\cite{Bardella2009, Bardella2010}, Gurtin \cite{GURT2004}, Neff et
al. \cite{NCA}, Ebobisse et al. \cite{EBONEFF, ENR2015, EHN2016});

\item[-] satisfy an extended positive dissipation principle and be therefore thermodynamically admissible;

\item[-]  be able to be cast in a convex analytical framework;

\item[-] support a well-posedness result in both the rate-independent
and the rate-dependent cases;

\item[-] have physical meaningful and transparent boundary conditions for the plastic varia\-bles.

\end{itemize} In this paper, we will propose such a model, which, in the end, has
a certain resemblance with the early model proposed by Menzel and
Steinmann \cite{MENZSTEIN} and we will show its well-posedness in
the rate-independent case with a suitable regularization. Our derivation of the model based on
invariance principles is, to our knowledge, entirely new. The model
is derived through a process of linearization of some state
variables in the finite strain case. The derivation of a linear
model from  a finite strain one is not new in the context of
elasto-plasticity. For instance, Mielke and Stefanelli
\cite{MielkeStefan2013} used $\Gamma$-convergence to derive
rigorously a model of linearized plasticity as limit of some finite
strain plasticity model. 

The mathematical well-posedness of our
model seems to be interesting in its own right. Once more, the
convex analytical framework, based on incorporating the postulate of
maximum plastic dissipation,\footnote{The postulate of maximum
plastic dissipation (PMPD), which was derived independently in
classical infinitesimal theory of plasticity from the so-called
Drucker's postulate by von Mises \cite{Mises1928}, Taylor
\cite{Taylor1947}, Hill \cite{Hill1948} Mandel \cite{Mandel1964}
(and later as a consequence of Il'iushin's postulate of plasticity
in strain space) has the form $\la\sigma
-\sigma^*,\dot{\bvarepsilon}_p\ra\geq0$, where $\sigma$ is the
actual stress tensor, $\dot{\bvarepsilon}_p$ is the plastic
strain-rate tensor and $\sigma^*$ is any admissible stress tensor.
The PMPD is equivalent to the associated flow rule in the dual
formulation in the local theory of plasticity. The maximal
dissipation (associated flow rule) simplifies the modelling
framework and facilitates the mathematical treatment.} put
forward initially by Moreau \cite{Moreau1976} and used later by many
authors (e.g.  \cite{Han-ReddyBook, DEMR1,
REM, NCA, EMR2008, EBONEFF, ENA2018, ENF2018, ENR2015, EHN2016}), proves to be ideally suited. For that purpose, our "additive"
model also admits a finite strain parent model by Krishnan and
Steigmann \cite{KRISSTEIG2014}, who did not consider, however, the
incorporation of (nonlocal) kinematic hardening.

Decisive for our new strain gradient plasticity model is the
introduction of Kr\"oner's incompatibility tensor $\inc$
(\cite{Kroner1955, Kroner1960, Kondo1955, Vangoethem2011,Vangoethem2016,Vangoethem2017}) as an
inhomogeneity measure acting on the symmetric plastic strain
$\bvarepsilon_p=\sym p$. This {\bf incompatibility tensor} is given by
\begin{equation}\label{inc-tensor1}
\inc(\sym p):=\Curl([\Curl \sym p]^T)\,\end{equation} and it
coincides to first order with the Riemann-Christoffel curvature
tensor $\mathcal{R}$ in the metric characterized by the finite plastic
strain tensor $C_p:={F_p\hskip-.1truecm}^T\,\hskip-.1truecm F_p$ (De
Wit \cite{Dewit1981}). In fact, considering the non-symmetric plastic distortion
 $F_p$ in the multiplicative decomposition (\ref{multi-split1}) and writing $F_p=\id +p$, the connection is
\begin{equation}\label{Riemtensor1}\mathcal{R}({F_p\hskip-.1truecm}^T\,\hskip-.1truecm F_p)=\mathcal{R}((\id
+p)^T(\id +p))=\mathcal{R}(\id +2\sym p+p^Tp)=2\inc(\sym
p)+\mbox{h.o.t.},\end{equation} where\footnote{Note that
$\mathcal{R}^i_{ jkl}$ is both geometrically and physically
nonlinear.}
\begin{equation}\label{def-Riem-curvature}
\mathcal{R}^i_{ jkl}(C_p):=\frac{\partial\Gamma^i_{jl}}{\partial
x_k}-\frac{\partial\Gamma^i_{jk}}{\partial x_l}
+\Gamma^m_{jl}\,\Gamma^i_{mk}-\Gamma^m_{jk}\,\Gamma^i_{ml}\,,
\end{equation}
with $$\Gamma^l_{ij}:=\frac12g^{kl}\left[\frac{\partial
g_{ki}}{\partial x_j}+\frac{\partial g_{jk}}{\partial
x_i}-\frac{\partial g_{ij}}{\partial x_k}\right]$$ being the
Christoffel symbols of the second kind in the metric $C_p=(g_{ij})$
whose inverse is $C_p^{-1}=(g^{ij})$. \vskip.2truecm \noindent Assume
that we change the reference configuration by a smooth invertible
map \begin{equation}\label{change-coord}\psi:\xi\mapsto
x(\xi)\,,\end{equation}
 then the plastic distortion $F_p$ in the
multiplicative decomposition (\ref{multi-split1}) should transform
according to \begin{equation}\label{transf-psip1}F_p(x)
\,\longrightarrow\, F_p(x(\xi))\,\frac{\partial x(\xi)}{\partial
\xi}=\widetilde{F}_p(\xi)\,\nabla_\xi\psi(\xi)\,.\end{equation}
 Based on the
Riemann-Christoffel curvature tensor $\mathcal{R}$ in
(\ref{def-Riem-curvature}), we may form a "true scalar" quantity,
namely the so-called {\it Lanczos scalar}
\begin{equation}\label{Lanczos1}\sum_{ijkl}\mathcal{R}^i_{ jkl}\cdot\mathcal{R}^i_{ jkl}=:
\norm{\mathcal{R}}^2_{\mathbb{R}^{81}}
\,,\end{equation} which is
 {\it form-invariant} under
a change of the reference configuration in the sense that
\begin{equation}\label{Riemtensor2}\norm{\mathcal{R}_x({F_p(x)}^T\hskip-.1truecm
F_p(x))}^2=\norm{\mathcal{R}_\xi((\widetilde{F}_p(\xi)\,\nabla_\xi\psi(\xi))^T(\widetilde{F}_p(\xi)\,\nabla_\xi\psi(\xi)))}^2\,,\end{equation}
(see Lanczos \cite[eq.\hskip.1truecm(2.3)]{Lanczos}) where $\mathcal{R}_x$ and
$\mathcal{R}_\xi$ denote the Riemann-Christoffel curvature tensor
expressed in coordinates $x$ and $\xi$, respectively. Moreover, if we consider the transformation
\begin{equation}\label{transf-direct-inv}F_p(x) \,\longrightarrow\,
\widehat{F}_p:=Q(x)F_p(x),\quad\mbox{with arbitrary }
Q(x)\in\Oo(3)\,,\end{equation} then we have the direct invariance
condition for the Riemann-Christoffel curvature tensor
\begin{eqnarray}\label{directinv1}\nonumber \mathcal{R}_x({F_p(x)}^T\hskip-.1truecm
F_p(x))&=&\mathcal{R}_x({F_p(x)}^T\hskip-.1truecm Q(x)^TQ(x)F_p(x))
=\mathcal{R}_x((Q(x)F_p(x))^T\hskip-.1truecm
Q(x)F_p(x))\\
&=&\mathcal{R}_x({\widehat{F}_p(x)}^T\hskip-.1truecm
\widehat{F}_p(x))\,.\end{eqnarray}

 \noindent The form-invariance property (\ref{Riemtensor2}) is
inherited in its geometric linearization, now by the $\inc$-operator, in the form of a direct invariance condition for the complete
operator (and not just the ``Lanczos''-type scalar $\norm{\inc_x\sym p}^2$):
\begin{equation}\label{inv-inc}\inc_x(\sym p(x))=\inc_x(\sym(p(x)
+\nabla_x\vartheta(x))\qquad\forall\,\vartheta(x)\in
C^\infty(\mathbb{R}^3,\,\mathbb{R}^3)\,,\end{equation} where we have
identified $\psi(x)=x+\vartheta(x)$.  
In fact, from the identity
(see \cite[Proposition 2.1]{GNMM2015})
\begin{equation}\label{inc-propperty1}
\nabla\axl(\skew\nabla\vartheta)=\frac12\nabla \curl
\vartheta=(\Curl \sym\nabla \vartheta)^T,\end{equation} taking
the Curl on both sides, we get
\begin{equation}\label{inc-property2}
0=\Curl\bigl(\frac12\nabla\curl \vartheta\bigr)=\Curl\bigl[(\Curl
\sym\nabla \vartheta)^T\bigr]=\inc(\sym\nabla
\vartheta)\,,\end{equation} which shows (\ref{inv-inc}). Of course
 (\ref{inv-inc}) implies that
\begin{equation}\label{inv-inc-lanczos}\norm{\inc_x(\sym p(x))}^2=\norm{\inc_x(\sym(p(x)
+\nabla_x\vartheta(x))}^2\qquad\forall\,\vartheta(x)\in
 C^\infty(\mathbb{R}^3,\,\mathbb{R}^3)\,,\end{equation}
 mirroring property (\ref{Riemtensor2}) for the Lanczos-type scalar.
 \vskip.2truecm
 \noindent In addition, the direct
invariance condition (\ref{directinv1}) for $\mathcal{R}$ under
rotation fields $Q(x)\in\SO(3)$ translates to an invariance
condition on the inc-operator as well. We write {$Q(x)=\id
+A(x)+\mbox{h.o.t.}$ with $A(x)\in\so(3)$} and in terms of the
infinitesimal plastic distortion, we consider
$$p(x)\,\longrightarrow\, A(x)+p(x)\qquad\forall A(x)\in\so(3)$$ and we
have the invariance condition
\begin{equation}\label{inc-inv1}\inc_x(\sym
p(x))=\inc_x(\sym(A(x)+p(x)))\qquad\forall\,A(x)\in\so(3)\,.\end{equation}
For both tensors $\mathcal{R}$ and $\inc$ we note the Saint-Venant
compatibility condition and its linearization:
\begin{equation}\label{transf-vanish-R-inc}\mbox{\begin{tabular}{|l|}\hline\\$\psi(x)=x+\vartheta(x)$\\\\
\hline\end{tabular}}  \quad\mbox{ and }\quad
   \begin{tabular}{|l|}\hline\\
 $C_p=F_p^T\hskip-.1truecm\, F_p=\id+\bvarepsilon_p+\mbox{h.o.t}$\\\\\hline
\end{tabular}   \end{equation}\vskip-.2truecm
$$\mbox{\begin{tabular}{|l|}\hline\\
$C_p=F_p^T\hskip-.1truecm\, F_p\in\mbox{Sym}^+(3)$ symmetric positive definite\\\\
$\mathcal{R}(C_p)=0\,\Leftrightarrow\,
C_p={\nabla\psi}^T\nabla\psi$\\\\
See e.g. Ciarlet-Laurent \cite[Theorem 1.1]{Ciarletlaurent2003}\\
\hline
\end{tabular}
$\approx$ \begin{tabular}{|l|}\hline\\
 $\inc(\sym p)=0\,\Leftrightarrow\,\sym p=\sym\nabla\vartheta$\\\\\hline
\end{tabular}}$$ in simply connected domains (see
\cite{MARSHUGH1983,Ciarlet2005,Ciarletlaurent2003,MAGSCAVANG}). For more properties of
the $\inc$-operator, we refer the reader to \cite{Vangoethem2011,
AMSVANGOETH, MAGSCAVANG}. \vskip.2truecm
  With these preliminaries, both tensors $\mathcal{R}$
and $\inc$ qualify as incompatibility measures on positive definite
symmetric plastic strains $C_p=F_p^T\hskip-.1truecm\, F_p$ in the
geometrically nonlinear and on symmetric plastic strains
$\bvarepsilon_p=\sym p$ in the geometrically linear settings,
respectively. \vskip.2truecm
 In a purely
phenomenological context, we have another way to measure the
incompatibility of the plastic distortion $F_p$ itself via the
so-called {\it dislocation density tensor} $\Curl_x F_p(x)$ (see for
instance \cite{GURTAN2005FD,NCA}).
 Following the works of Davini-Parry \cite{DAVPAR}, Cermelli-Gurtin
\cite{CERGURT2001} and Epstein \cite{EPSTEIN2010}, it has been shown
that the differential operator (the "true dislocation density
tensor")\,\footnote{The role of that tensor has been critically
discussed in Acharya \cite{Acharya2008}.}
\begin{equation}\label{true-disl-density}
\frac1{\det F_p}\,(\Curl F_p)\,F^T_p=\det
F_e(\Curl_CF_e^{-1})F_e^{-T}\,,\end{equation} (where $\Curl_C$ means
taking the spatial Curl w.r.t. the current configuration) is
 also {\it form-invariant} under arbitrary coordinate transformation of the
reference placement (see \cite[eq.(1.4) on
p.1542]{CERGURT2001}).\footnote{Note that Gurtin uses a different
definition of the Curl-operator. We have the relation $\Curl
F_p=[\Curl_{\mbox{\scriptsize{Gurtin}}}F_p]^T$ (see
\cite{GURTAN-BOOK}).} This property is easily understood from
(\ref{true-disl-density})$_2$, which is invariant under
reparame\-trization of the reference system anyway.

 \noindent From now on, we assume plastic
incompressibility, i.e., $\det F_p=1$ and $\tr p=0$.
Either\begin{equation}\label{curl-incomp1} (\Curl_xF_p)\,F^T_p
\quad\mbox{linearized through}\quad F_p=\id+p \quad\mbox{turns
into}\quad \Curl_xp\end{equation} or
\begin{equation}\label{inc-incomp2}\mathcal{R}({F_p\hskip-.1truecm}^T\hskip-.1truecm
F_p)\quad\mbox{linearized through}\quad F_p=\id+p \quad\mbox{turns
into}\quad\inc(\sym p)\,.\end{equation} Therefore, in the
geometrically linear setting two measures of incompatibility of the
infinitesimal plastic distortion $p$ are used in the literature.
Some authors use the tensor $\Curl p$ (see \cite{berdisedov1967,
Berdi2006, GURTAN2005, NCA}) while others use $\inc(\sym p)$ (see
\cite{MENZSTEIN}).

 \begin{table}[h!]\small
\begin{center}\begin{tabular}{|l|l|}\hline &\qquad\\
$\begin{array}{lcl}\Curl p=0 &\Rightarrow& p=\nabla\vartheta\\
 &\Rightarrow& \inc(\sym p)=\inc(\sym\nabla\vartheta)=0\end{array}$ & 
$\begin{array}{lcl}\inc(\sym p)=0 &\Rightarrow& \sym p=\sym\nabla\vartheta\\
 &/\hskip-.35truecm\Rightarrow&  p=\nabla\vartheta\,\,\mbox{ but}\\
 & \Rightarrow & p=\nabla\vartheta +A(x),\quad A(x)\in\so(3)\\
  &/\hskip-.35truecm\Rightarrow&\Curl p=0\end{array}$ 

\\
 \hline
\end{tabular}\end{center}\caption{\footnotesize  We need to realize that $\Curl p$ is a ``sharper'' incompatibility measure than $\inc(\sym p)$. }\label{table:curl-inc}\end{table}

 Concerning the invariance of the curvature tensor $\mathcal{R}$ observed in
(\ref{Riemtensor2}), we note that under a change of reference
placement
$F_p(x)\,\longrightarrow\,\widetilde{F}_p(\xi)\,\nabla\psi(\xi)$,
we obtain directly the {\it form-invariance} of the true dislocation density
tensor $(\Curl F_p)\,F^T_p$ as well, meaning that
\begin{equation}\label{form-inv1-disloc1}
(\Curl_xF_p)\,F^T_p=(\Curl_\xi(\widetilde{F}_p(\xi)\,\nabla\psi(\xi)))\,(\widetilde{F}_p(\xi)\,\nabla\psi(\xi))^T\,,
\end{equation}
 where $\Curl_x$ and $\Curl_\xi$ are the Curl expressed in
coordinates $x$ and $\xi$, respectively. Therefore, the expression
$\norm{(\Curl_xF_p)\,F^T_p}^2$ is also  a true scalar quantity. Accordingly, in
the linearized setting we consider as in (\ref{inv-inc})
\begin{equation}\label{change-ref-place-lin1}p(x)\,\longrightarrow\, p(x) +\nabla\vartheta(x)\end{equation} and we
obtain directly  the invariance
\begin{equation}\label{form-inv1-disloc-lin1}
\Curl p(x)=\Curl \bigl[p(x) + \nabla\vartheta(x)\bigr]\,,\end{equation}
 similar to (\ref{inv-inc}).

\noindent Concerning superposition of rotation fields, we note that
for $F_p\,\longrightarrow\, \overline{Q}\,F_p$ for
$\overline{Q}\in\Oo(3)$ where $\overline{Q}$ is a homogeneous
rotation, we have
\begin{equation}\label{rotation-inv1-disloc1}
(\Curl(\overline{Q}\,F_p))(\overline{Q}\,F_p)^T=\overline{Q}\,(\Curl
F_p)F^T_p\,{\overline{Q}}^T\,.\end{equation} Similarly, in the
linearized setting we consider $p(x)\,\longrightarrow\,
\overline{A}+p(x)$ for every $\overline{A}\in\so(3)$ with
$\overline{A}$ a constant skew-symmetric matrix and we get
\begin{equation}\label{inv-superposition-skew}\Curl p(x)=\Curl\bigl[\overline{A}+p(x)\bigr]\,.\end{equation} We
have the relation
$$\mbox{\begin{tabular}{|l|}\hline\\$\psi(x)=x+\vartheta(x)$\\\\
\hline\end{tabular}}\quad\mbox{ and }\quad
\begin{tabular}{|l|}\hline\\
$F_p=\id+p$\\\\
\hline
\end{tabular}$$\vskip.2truecm
\begin{equation}\label{vanish-curl}\mbox{\begin{tabular}{|l|}\hline\\
$(\Curl F_p)F^T_p=0\,\Leftrightarrow\,
F_p=\nabla\psi$\\\\
\hline
\end{tabular}\quad $\equiv$\quad
\begin{tabular}{|l|}\hline\\
$\Curl p=0\,\Leftrightarrow\, p=\nabla\vartheta$\\\\
\hline
\end{tabular}}\end{equation} in simply connected domains and under appropriate regularity conditions  (see
\cite[Section 59]{HAY1953}).
\vskip.2truecm
 The difference between the two incompatibility measures for the linearized setting,
$\inc(\sym p)$ on the one hand, and $\Curl p$ on the other hand is
the invariance property under superposed infinitesimal rotations:
 while $\inc(\sym p)$ allows to {\bf superpose any inhomogeneous}
infinitesimal  rotation field $A(x)$, $\Curl p$ allows to
{\bf superpose only homogeneous} infinitesimal rotation fields  $\overline{A}$.

In single crystal gradient plasticity, it is typically $\Curl
\hskip-.06truecm F_p$ which is used whereas it is debatable, whether
$\Curl\hskip-.06truecm F_p$ is a good state-variable for
polycrystalline material without texture.
\vskip.2truecm
In the following we will use a set of invariance conditions which
will allow us to decide between using $\Curl$ or $\inc$. It is
clear, however, that assuming a set of invariance requirements is
already a constitutive requirement and therefore subject to
discussion.{\footnote{Only time-objectivity requirements are not to
be discussed.}} In all these developments, beyond the discussion on
which invariance principles are applied, it is our aim to clearly
state and show, which kind of modelling restrictions will be
obtained from them}.
 \vskip.2truecm  Our contribution is structured as follows: after introducing in Section \ref{Notations}  some notations, operators and function spaces used throughout the paper, we set the stage in Section \ref{inv-conditions} with two important invariance conditions on which our  model will be tested. Namely, the {\it gauge-invariance} known as invariance under compatible transformations of the reference system, and a  novel rotational invariance
postulate for polycrystals, called {\it micro-randomness}. In Section \ref{some models},  we first present few models of gradient plasticity with Kr\"{o}ner's
incompatibility tensor which fail our invariance conditions, then we introduce our novel fourth order phenomenological model which, though it fails also the gauge-invariance condition, is invariant w.r.t. a subclass of reparametrizations of the reference configurations, including the infinitesimal conformal group. The new model is then formulated using the convex analytical framework leading to  mathematical strong and weak formulations. Finally an  existence result for the weak formulation is obtained.
\vskip.2truecm\noindent
Let us next
fix some notations and definitions which will also make the paper
more clear and readable.
\section{Some notational agreements and
definitions}\label{Notations} Let $\Omega$ be a bounded domain
 in $\BBR^3$ with Lipschitz continuous boundary $\partial\Omega$, which is occupied by an elastoplastic
body in its undeformed configuration. Let $\Gamma$ be a smooth
subset of $\partial\Omega$ with non-vanishing $2$-dimensional
Hausdorff measure. A material point in $\Omega$ is denoted by $x$
and the time domain under consideration is the interval $[0,T]$.\\
 For every $a,\,b\in\BBR^3$, we let $\la a,b\ra_{\BBR^3}$ denote the scalar
 product on $\BBR^3$ with associated vector
norm $\norm{a}^2_{\BBR^3} = \la a, a\ra_{\BBR^3}$. We denote by
$\BBR^{3\times 3}$ the set of real $3\times 3$ tensors. The standard
Euclidean scalar product on $\BBR^{3\times 3}$ is given by $\la
A,B\ra_{\BBR^{3\times 3}} = \mbox{tr}\,\bigl[AB^T\bigr]$, where
$B^T$ denotes the transpose tensor of $B$. Thus, the Frobenius
tensor norm is $\norm{A}^2 = \la A,\,A\ra_{\BBR^{3\times 3}}$.
 In the following we omit the subscripts $\BBR^3$ and $\BBR^{3\times 3}$. The identity tensor on $\BBR^{3\times 3}$ will be denoted by
  $\id$, so that $\mbox{tr}(A) = \la A, \id\ra$. We let $\mbox{GL}(3):=\{X\in\mathbb{R}^{3\times
3}\,\,|\,\,\det(X)\neq0\}$ denote the group of invertible $3\times
3$ square matrices; $\mbox{GL}^+(3):=\{X\in\mathbb{R}^{3\times
3}\,\,|\,\,\det(X)>0\}$; $\SO(3):=\{X\in
\mbox{GL}(3)\,\,|\,\,X\,X^T=\id\,,\,\,\det[X]=1\}$ is the Lie Group of
rotations in $\mathbb{R}^3$ whose Lie Algebra is the set
$\so(3):=\{X\in\BBR^{3\times 3}\,|\,\,X^T=-X\}$ of skew-symmetric
tensors.
 We let
$\mbox{Sym\,}(3):=\{X\in\BBR^{3\times 3}\,|\,\,X^T=X\}$ denote the
set of symmetric tensors and $\sL(3):=\{X\in\BBR^{3\times
3}\,|\,\,\mbox{tr\,}(X)=0\}$ be the Lie Algebra of traceless
tensors. For every $X\in\BBR^{3\times 3}$, we set
$\sym(X)=\frac12\bigl(X+X^T\bigr)$,
$\skew\,(X)=\frac12\bigl(X-X^T\bigr)$ and
$\dev(X)=X-\frac13\mbox{tr}\,(X)\,\id\in\sL(3)\,$ for the symmetric
part, the skew-symmetric part and the deviatoric part of $X$,
respectively. Quantities which are constant in space will be denoted
with an overbar, e.g., $\overline{A}\in\so(3)$ for the function
$A:\mathbb{R}^3\to\so(3)$ which is constant with constant
$\overline{A}$.

The body is assumed to undergo deformations. Its behaviour is
governed by a set of equations and constitutive relations. Below is
a list of variables and parameters used throughout the
paper:\begin{itemize}
\item[$\bullet$] $\varphi$ is the deformation of the body;
\item[$\bullet$] $u(x,t)=\varphi(x,t)-x$  is the displacement of the macroscopic material
points;
\item[$\bullet$] $F=\nabla \varphi=\id+\nabla u$ is the deformation
gradient;
\item[$\bullet$] $F_p=\id+p$ is the plastic distortion which is a non-symmetric tensor with unit determinant, that is, $F_p\in\SL(3)$;

\item[$\bullet$] $F_e=\id +e$ is the elastic distortion which  is a non-symmetric tensor;
\item[$\bullet$] $C_p:=F_p^T\,F_p=\id+\bvarepsilon_p+\ldots\,$ is the positive definite plastic metric;
\item[$\bullet$] $C_e:=F_e^T\,F_e=\id+\bvarepsilon_e+\ldots\,$ is the positive definite elastic strain tensor;
\item[$\bullet$] $p$ is the infinitesimal plastic distortion variable which  is a
non-symmetric second order tensor, incapable of sustaining
volumetric changes; that is, $p\in\sL(3)$. The tensor $p$\,
represents the average plastic slip; $p$ is not gauge-invariant, while the rate $\dot{p}$ is;

\item[$\bullet$] $e=\nabla u -p$ is  the infinitesimal elastic distortion which  is a
non-symmetric second order tensor and is a state-variable;

\item[$\bullet$] $\bvarepsilon_p=\sym p$ is the symmetric infinitesimal plastic strain
tensor, which is also trace free, {$\bvarepsilon_p\in\sL(3)$;} $\bvarepsilon_p$ is not gauge-invariant; the rate $\dot{\bvarepsilon}_p=\sym\dot{p}$ is gauge-invariant; $\bvarepsilon_p$ is not a state-variable;

\item $\skew p$ is called plastic rotation or plastic spin and is not a state-variable;

\item[$\bullet$] $\bvarepsilon_e=\sym\,(\nabla u -p)$ is the symmetric infinitesimal  elastic
strain tensor and is a state-variable;

\item[$\bullet$] $\sigma$ is the Cauchy stress tensor which  is a symmetric
second order tensor and is gauge-invariant; 

\item[$\bullet$] $\yieldzero$ is the initial
yield stress for plastic strain  and is gauge-invariant;

\item[$\bullet$]  $\yieldlimit$ is the current yield stress for plastic strain  and 
 is gauge-invariant;

\item[$\bullet$] $f$ is the body force;

\item[$\bullet$] $\Curl p=-\Curl e=\alpha$ is Nye's dislocation density
tensor (see (\ref{true-disl-density}) for the definition of
$\alpha$), satisfying the so-called Bianchi identities
$\Div\alpha=0$ and is gauge-invariant;

\item[$\bullet$] $\mathcal{R}$ is the Riemann-Christoffel curvature tensor, see (\ref{def-Riem-curvature});

\item[$\bullet$] $\norm{\mathcal{R}}^2$ is the Lanczos-type scalar;

\item[$\bullet$] $\inc(\sym p)=\inc(\bvarepsilon_p)=-\inc(\bvarepsilon_e)$ is
Kr\"{o}ner's second order incompatibility tensor and is gauge-invariant;

\item[$\bullet$] $\gamma_p=\dsize\int_0^t\norm{\sym\dot{p}}\,ds$ is the accumulated equivalent  plastic strain and is gauge-invariant;
\item[$\bullet$] $\widetilde{\gamma}_p=\dsize\int_0^t\norm{\dot{p}}\,ds$ is the accumulated equivalent  plastic distortion and is gauge-invariant.
\end{itemize}
\vskip.2truecm\noindent For isotropic media, the fourth order
isotropic elasticity tensor $\C_{\mbox{\footnotesize{iso}}}:\mbox{Sym}(3)\to\mbox{Sym}(3)$ is
given by
\begin{equation}
\C_{\mbox{\footnotesize{iso}}}X = 2\mu\,\dev\,\sym X+\kappa
\,\tr(X) \id =2\mu\,\sym X+\lambda\,\tr(X)\id\label{C}
\end{equation}
for any second-order tensor $X$, where $\mu$ and $\lambda$ are the
Lam{\'e} moduli satisfying
\begin{equation}\label{Lame-moduli}
\mu>0\quad\mbox{ and }\quad 3\lambda +2\mu>0\,,
\end{equation} and $\kappa>0$ is the bulk modulus.
These conditions suffice for pointwise ellipticity of the elasticity
tensor in the sense that there exists a constant $m_0 > 0$ such that
\begin{equation}
\forall X\in\mathbb{R}^{3\times 3}\mbox{: }\quad\la\sym X,\C_{\mbox{\footnotesize{iso}}}\sym X\ra \geq m_0\, \norm{\sym X}^2\,.
\label{ellipticityC}
\end{equation}

For every $X\in C^1(\Omega,\,\BBR^{3\times 3})$ with rows
$X_1,\,X_2,\,X_3$, we use in this paper the definition of $\Curl X$
in \cite{NCA, SVEN}:
\begin{equation}\label{def-Curl}\Curl X =\left(\begin{array}{l}\mbox{curl\,}X_1\,\,-\,\,-\\
\mbox{curl\,}X_2\,\,-\,\,-\\
\mbox{curl\,}X_3\,\,-\,\,-\end{array}\right)\in\BBR^{3\times
3}\,,\end{equation} for which $\Curl\,\nabla v=0$ for every $v\in
C^2(\Omega,\,\BBR^3)$. Notice that the definition of $\Curl\,X$
above is such that $(\Curl X)^Ta=\mbox{curl\,}(X^Ta)$ for every
$a\in\BBR^3$ and this clearly corresponds to the transpose of the
Curl of a tensor as defined in
\cite{GURTAN2005, GURTAN-BOOK}.\\

 For
$$
\overline{A}=\left(\begin{array}{ccc}
0 &-a_3&a_2\\
a_3&0& -a_1\\
-a_2& a_1&0
\end{array}\right)\in \so(3)\,,
$$
we consider the operator $\axl:\so(3)\rightarrow\mathbb{R}^3$ and
$\anti:\mathbb{R}^3\rightarrow \so(3)$ through

\begin{equation}\label{axl-anti1}
\axl(\overline{A}):=\left( a_1, a_2, a_3 \right)^T,\quad \quad
\overline{A}.\, v=(\axl \overline{A})\times v, \quad
(\anti(v))_{ij}=\varepsilon_{jik}\,v_k, \quad \forall \,
v\in\mathbb{R}^3,\end{equation}\begin{eqnarray}\label{axl-anti2}
(\axl \overline{A})_k &=&-\frac{1}{2}\,\sum\limits_{i,j=1}^3
\epsilon_{ijk}\,\overline{A}_{ij}=\frac{1}{2}\,\sum\limits_{i,j=1}^3
\epsilon_{kij}\,\overline{A}_{ji}\,,\\
\overline{A}_{ij}&=&\sum\limits_{k=1}^3-\epsilon_{ijk}\,(\axl
\overline{A})_k=:\anti(\axl \overline{A})_{ij}\,,\\
\label{anti-axl1}\axl(\anti(v))_k&=&v_k\,,
\end{eqnarray}
 where $\epsilon_{ijk}$ is the totally antisymmetric third order Levi-Civita permutation
 tensor defined by
$$\epsilon_{ijk}:=\left\{\begin{array}{ll} 1 &\mbox{ if
}\{i,j,k\}=\{1,2,3\}\,,\{2,3,1\}\mbox{ or }\{3,1,2\}\,,\\
-1 &\mbox{ if }\{i,j,k\}=\{2,1,3\}\,,\{1,3,2\}\mbox{ or }\{3,2,1\}\,,\\
0 &\mbox{ if an index is repeated}.\end{array}\right.$$ Hence, the
operators $\axl:\so(3)\to\mathbb{R}^3$ and
$\anti:\mathbb{R}^3\to\so(3)$ are canonical identifications of
$\so(3)$ and $\mathbb{R}^3$. Notice that,
\begin{eqnarray}\label{prop-axl1}\nonumber (\axl\skew A)_k=\frac{1}{2}\,\sum\limits_{i,j=1}^3
\epsilon_{kij}\skew(A)_{ji}&=&\frac{1}{4}\,\sum\limits_{i,j=1}^3
\epsilon_{kij}A_{ji}-\frac{1}{4}\,\sum\limits_{i,j=1}^3
\epsilon_{kij}A_{ji}\\
&=&\frac{1}{2}\,\sum\limits_{i,j=1}^3
\epsilon_{kij}A_{ji}\,\qquad\forall\,A\in\BBR^{3\times
3}\,.\end{eqnarray}

The following function spaces and norms will also be used later.
\begin{eqnarray}\label{Curl-spaces}
\nonumber \mbox{H}(\mbox{Curl};\,\Omega,\,\BBR^{3\times 3})&:=&\Bigl\{X\in
L^2(\Omega,\,\BBR^{3\times 3})\,\,|\,\,\mbox{Curl\,}X\in
L^2(\Omega,\,\BBR^{3\times 3})\Bigr\}\,,\\
\norm{X}^2_{\mbox{\scriptsize H}(\mbox{\scriptsize
Curl};\Omega)}&:=&\norm{X}^2_{L^2(\Omega)}+\norm{\mbox{Curl\,}X}^2_{L^2(\Omega)}\,,\quad\forall
X\in\mbox{H}(\mbox{Curl;\,}\Omega,\,\BBR^{3\times 3})\,,\\
\nonumber
\mbox{H}(\mbox{Curl};\,\Omega,\,\mathbb{E})&:=&\Bigl\{X:\Omega\to\mathbb{E}\,\,|\,\,X\in
\mbox{H}(\mbox{Curl};\,\Omega,\,\BBR^{3\times 3})\Bigr\}\,,
\end{eqnarray}
for $\mathbb{E}:=\sL(3)$ or
$\mbox{Sym}\,(3)\cap\sL(3)$.\vskip.2truecm\noindent
 We also consider the space
\begin{equation}\label{null-trace-space}\mbox{H}_0(\mbox{Curl};\,\Omega,\,\Gamma,\BBR^{3\times 3})\end{equation} as the completion in
the norm in  (\ref{Curl-spaces}) of the space $\{q\in
C^\infty(\Omega,\,\Gamma,\,\BBR^{3\times
3})\,|\,\,q\times n|_\Gamma=0\}\,.$ Therefore, this space
generalizes the tangential Dirichlet boundary condition
\begin{equation}\label{tang-boundary}q\times n|_\Gamma=0\,\end{equation}
to be satisfied by the plastic distortion $p$ or the plastic strain
$\bvarepsilon_p:=\sym p$. Whenever, $\Gamma=\partial\Omega$, we simply write $\mbox{H}_0(\mbox{Curl};\,\Omega,\,\BBR^{3\times 3})$. The space
$$\mbox{H}_0(\mbox{Curl};\,\Omega,\,\Gamma,\mathbb{E})$$ is defined
as
 in (\ref{Curl-spaces}). \vskip.2truecm\noindent
 The divergence operator Div on second order
tensor-valued functions is also defined row-wise as
\begin{equation}\label{def-div}\mbox{Div}\,X=\left(\begin{array}{l}\mbox{div\,}X_1\\
\mbox{div\,}X_2\\
\mbox{div\,}X_3\end{array}\right)\,.\end{equation} 
Further properties of Kr\"{o}ner's incompatibility tensor inc can be found in the appendix.
\section{Discussion of some invariance conditions in plasticity}\label{inv-conditions}
\subsection{Gauge-invariance - invariance under compatible transformations of the reference system}\label{gauge-inv}
Since the modelling should be invariant with respect to
the used coordinates system, we may introduce the {\it
gauge-invariance} condition.\\
Consider again the multiplicative split $F(x)=F_e(x)F_p(x)$ and
perform a compatible change of the reference configuration, i.e.,
set $x=\psi(\xi)$. Then we have upon transforming to new
coordinates\begin{equation}\label{gauge-inv0}
F(\psi(\xi))\,\nabla_\xi\psi(\xi)=F_e(\psi(\xi))\,F_p(\psi(\xi))\,\nabla_\xi\psi(\xi)\,.\end{equation}
Therefore we require our new model to be {\bf form-invariant} under
$$\boxed{\begin{array}{llll}
 F(x)& \longrightarrow & F(\psi(\xi))\nabla_\xi\psi(\xi)& \forall\psi\in C^\infty(\mathbb{R}^3,\mathbb{R}^3)\\ &&&\\
F_p(x) &\longrightarrow &
F_p(\psi(\xi))\nabla_\xi\psi(\xi)&\mbox{({\bf F}inite {\bf G}auge-{\bf I}nvariance)}\end{array}}\,.\eqno\mbox{
({\bf FGI})}
$$

Performing a geometrical linearization, we obtain
$$F(x)=F_e(x)F_p(x)\quad\longrightarrow\quad \nabla
u=e(x)+p(x)$$ and the finite gauge-invariance ({\bf FGI}) translates into
direct invariance under $$\boxed{\begin{array}{llll}
 \nabla u(x) & \longrightarrow & \nabla u(\xi) +\nabla\vartheta(\xi)&\forall\vartheta\in C^\infty(\mathbb{R}^3,\mathbb{R}^3)
\\&&&\\
p(x)&\longrightarrow & p(\xi)+\nabla\vartheta(\xi)
&\mbox{({\bf L}inear {\bf G}auge-{\bf I}nvariance)}\end{array}}\,,\eqno\mbox{({\bf LGI})}$$  which is also known
as {\bf translational gauge-invariance} (see Lazar \cite{LAZAR2000,
LAZAR2002, LAZANAS2008}).
\subsection{Micro-randomness: a novel rotational invariance
postulate for polycrystals}\label{micro-randomness} Polycrystals can
be viewed as random aggregates of single crystals  which, at
sufficiently large scales can be viewed as isotropic.

\begin{figure}[h!]
\centering
\vspace{-2.5cm}
\hspace{-1.7cm} \includegraphics[width=0.6\textwidth]{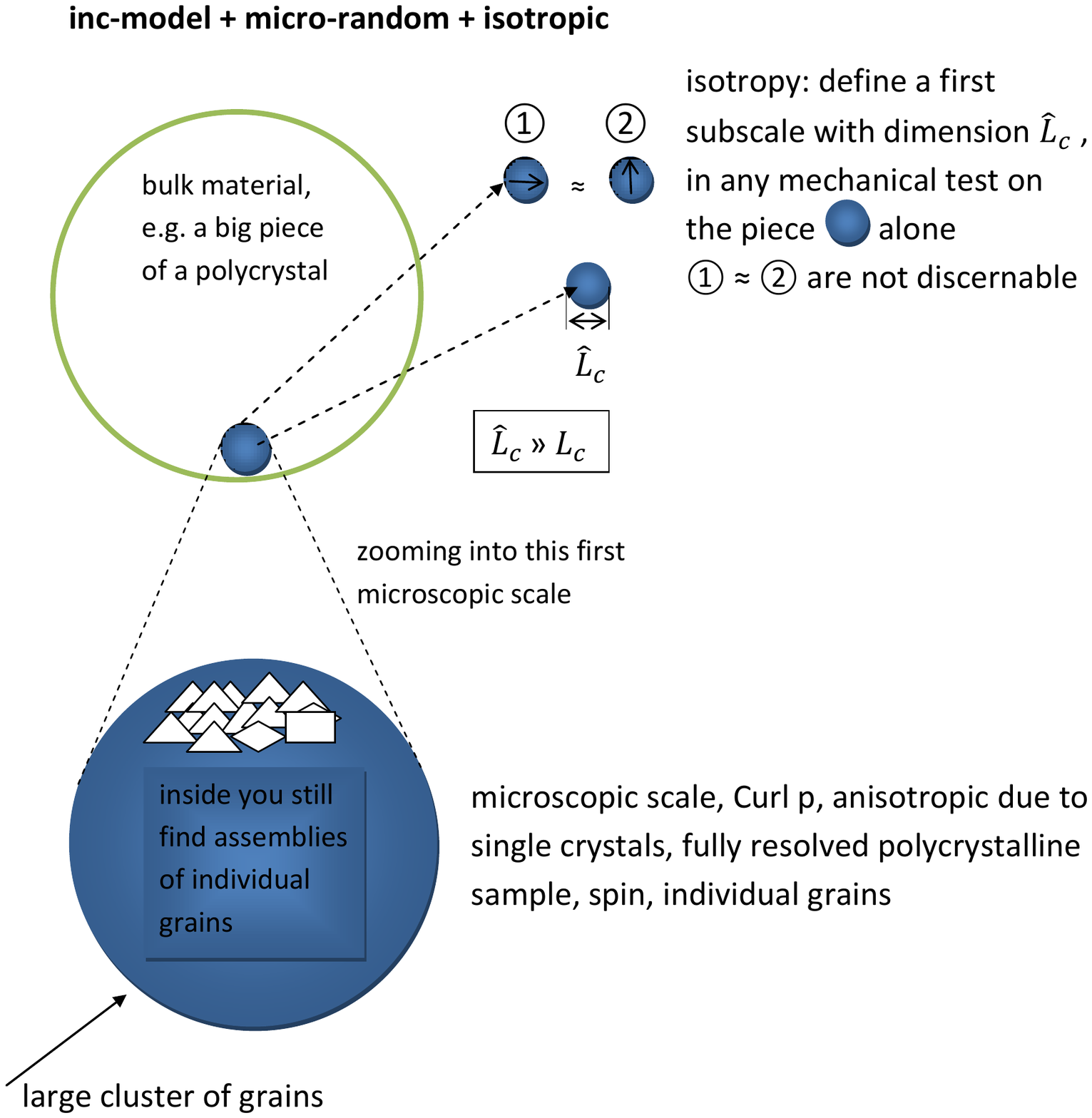}\hspace{-.5cm}       
\includegraphics[width=0.5\textwidth,   height=0.58\textheight]{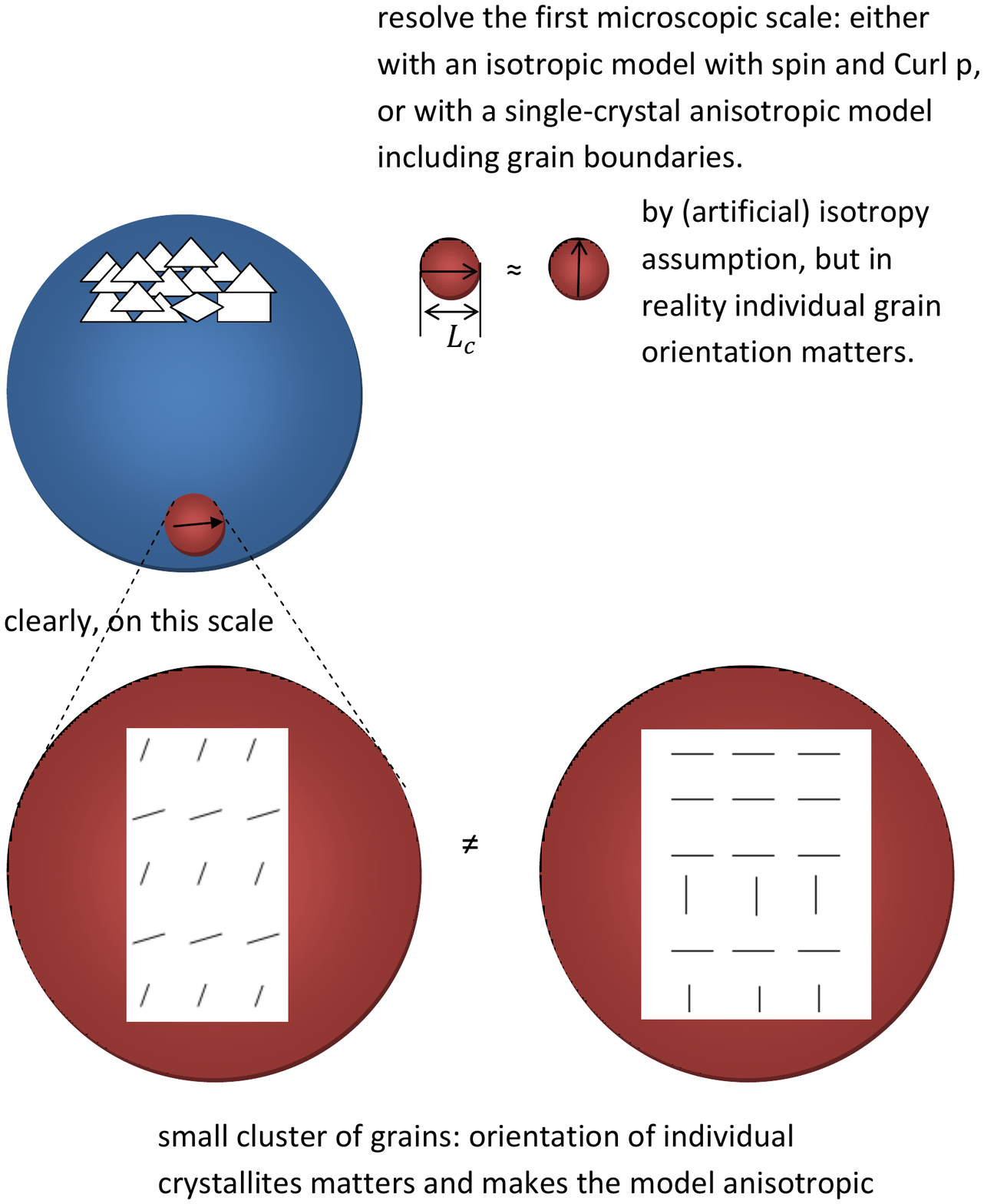}\vspace{-3.5cm}
 \includegraphics[width=0.5\textwidth,  height=0.6\textheight]{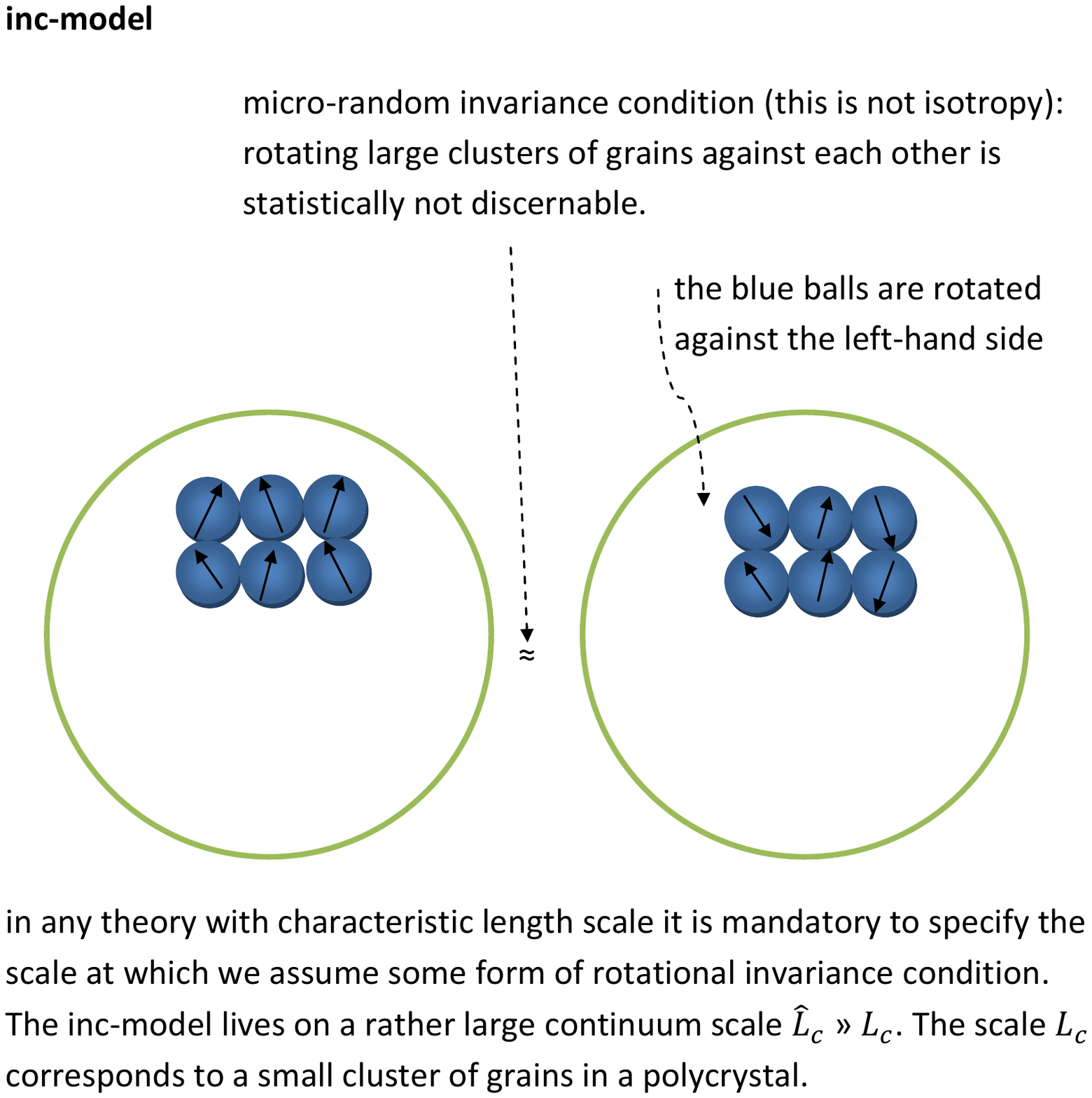}\hspace{-.5cm}
 \includegraphics[width=0.5\textwidth,  height=0.52\textheight]{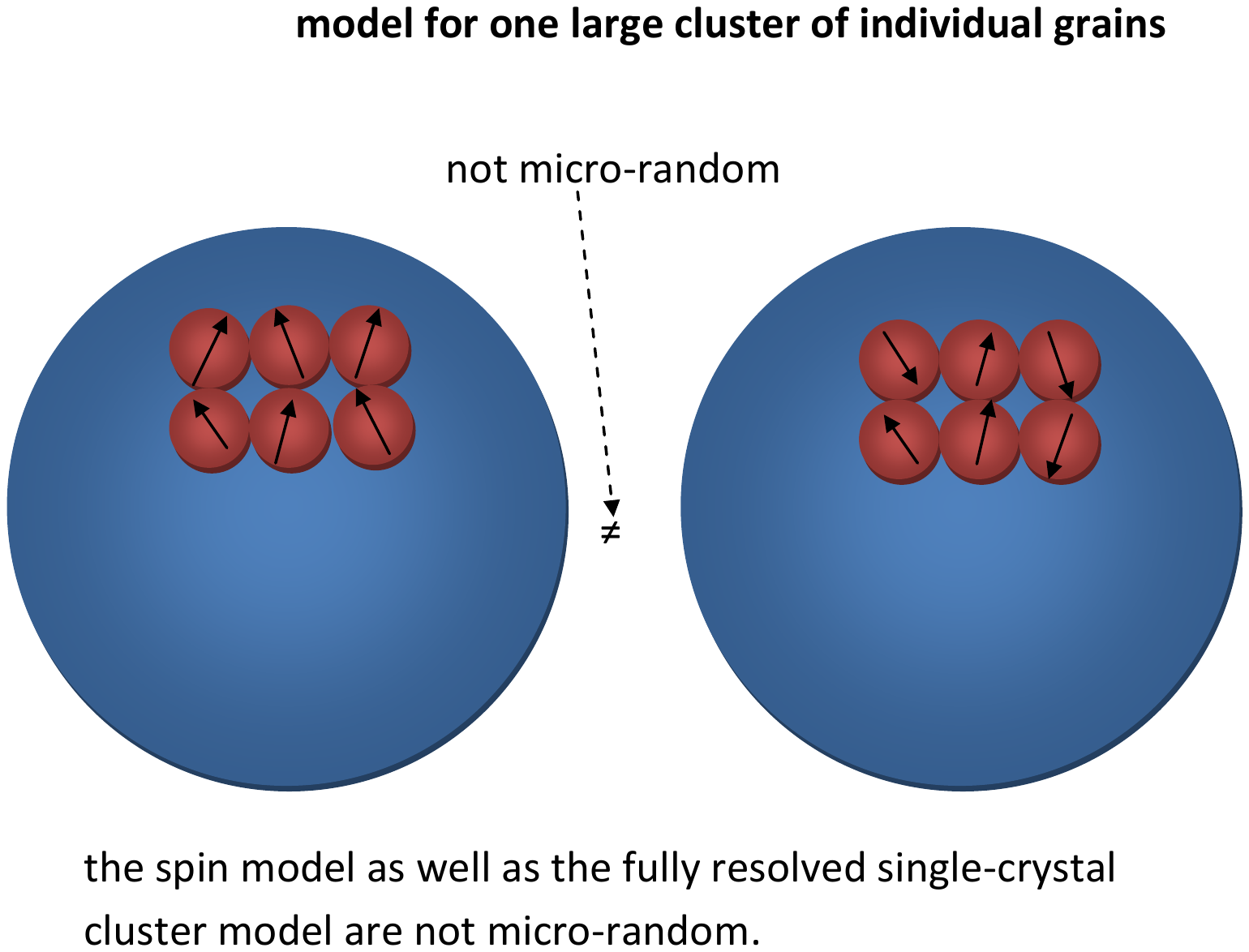}
\vspace{-4.2cm}
\caption{\footnotesize The micro-random model does not resolve a scale smaller than $\widehat{L}_c$ (blue-filled balls)! Inside the blue-filled balls, we have crystallites; if we want to model what happens inside there, we can still try to model this with isotropy assumption. This is our isotropic $\mbox{Curl}\, p$-model with spin (see Table \ref{table:isohard-fullyiso-new}). In reality however, the orientation of the crystallites matters.}
\label{fig:inc-microrandom}
\end{figure}

Imagine a
given initial distribution $F_p(x,0)$ of grains and subject the
polycrystal to a given mechanical loading which alters the plastic
state. The result will be recorded in the history $t\to F_p(x,t)$.

Now consider a randomly rotated initial distribution of grains and
plastic distortions via
\begin{equation}\label{initial-change}\boxed{
\widetilde{F}_p(x,0):=Q(x)F_p(x,0)\quad\mbox{with}\quad
Q(x)\in\SO(3)}\,.\end{equation}
 At a sufficiently large scale we are
not able to discern this rotational rearrangement and we are led to
assume that the new plastic history under the same given loading as
before should be
\begin{equation}\label{instant-t-change}\boxed{
t\to\widetilde{F}_p(x,t):=Q(x)F_p(x,t)}\,.\end{equation}

This is essentially a new invariance requirement to be imposed on
our model for the polycrystal. It means that, up to the initially
different inhomogeneous rotation of the grains, the response is the
same. Since rotations are involved one might take this as a
statement of classical isotropy. However, this would be misguided
since classical isotropy is concerned with rigidly rotating the
whole (polycrystalline) sample, while here each individual grain is
rotated differently\footnote{Rotating grains against each other
(see Neff et al. \cite{NJR2009}) in a polycrystal changes the
eigen-stresses along grain boundaries. Therefore, our new invariance
requirement cannot be a fundamental law of nature but may rather
serve to concentrate on some effective macroscopic features in a
homogenized model.} (see explanation in  Figure \ref{fig:inc-microrandom}).

Considering now the geometrically linear setting, we compare the
initial infinitesimal plastic distortion $p(x,0)$ with solution $t\to p(x,t)$
versus $p(x,0)+A(x)$ with its time evolution
$t\to\widetilde{p}(x,t)$. Our {\it micro-randomness invariance}
condition postulates in the geometrically linear setting that
$$\boxed{\begin{array}{lllllll}p(x,0)&\to&\widetilde{p}(x,0)&=&p(x,0)&+&A(x)\quad\Rightarrow\quad\\
p(x,t)&\to&\widetilde{p}(x,t)&=&p(x,t)&+&A(x)\quad\forall
t>0\end{array} \quad\mbox{({\bf L}inear {\bf M}icro-{\bf
R}andom)}}\eqno\mbox{\scriptsize\bf(LMR)}$$ and in the finite
deformation setting that
$$\boxed{\begin{array}{lll}F_p(x,0)&\to&\widetilde{F}_p(x,0)=Q(x)\,F_p(x,0)\quad\Rightarrow\quad\\
F_p(x,t)&\to&\widetilde{F}_p(x,t)=Q(x)\,F_p(x,t)\quad\forall
t>0\end{array} \quad\mbox{({\bf F}inite {\bf M}icro-{\bf
R}andom)}}\eqno\mbox{\scriptsize\bf(FMR)}$$
 In addition, we are aware of the fact that determining exact initial
conditions for the rotations of grains in a polycrystal is
practically impossible. Therefore, the influence of considering
different initial grain distributions should be minimized in
order to obtain a suitable effective model. Our micro-randomness
invariance condition ensures that the effect of different initial
conditions shows only as an "offset" of an otherwise unique response,
as seen above.\vskip.2truecm\noindent
 Notice that
\begin{equation}\label{micro-random-iso-finite}\left.\begin{array}{lclcl} F(x)&=&F_e(x)\,F_p(x) &=&
F_e(x)\underbrace{Q(x)^TQ(x)}_{=\,\id}F_p(x)\\&&&&\\
\nabla u&=&e(x)+p(x) &=& e+\underbrace{A^T(x)
+A(x)}_{=\,0}+p(x)\end{array}\right\}\,\qquad Q(x)=\exp(A(x))\,,\end{equation} thus the
invariance condition connected to {\bf micro-randomness} reads in the
finite strain case
\begin{equation}\label{micro-random-finite}(F,F_e,F_p)\quad\longrightarrow\quad (F,F_eQ^T(x),QF_p(x))\quad\forall\,
Q(x)\in \SO(3)\,,\end{equation} while in the geometrically linear
context, we need to require the invariance
\begin{equation}\label{micro-random-linear}(\nabla
u,e,p)\quad\longrightarrow\quad (\nabla u, e+A^T(x),
A(x)+p(x))\quad\forall \,A(x)\in \so(3)\,.\end{equation}
\subsection{Isotropy in geometrically linear models}\label{Isotropy} While the above invariance conditions can be characterized by additive operators, for classical  isotropy we need the group of rotations $Q\in\Oo(3)$. We define isotropy  in geometrically linear models to be {\bf form-invariance} under simultaneous change of spatial and referential coordinates by a rigid rotation. In this case, scalar functions $ h:\mathbb{R}^3\to\mathbb{R}$, vector fields $\phi:\mathbb{R}^3\to\mathbb{R}^3$ and second order tensor fields  $S:\mathbb{R}^3\to\mathbb{R}^{3\times 3}$ are transformed as follows:
\begin{equation}\label{sharp-transfo}\left.\begin{array}{lclcl} h &\longrightarrow & h^\sharp &\mbox{with}& h^\sharp(\xi):=h(Q^T\xi)\,,\\&&&&\\
\phi&\longrightarrow & \phi^\sharp &\mbox{with}& \phi^\sharp(\xi):=Q\,\phi(Q^T\xi)\,,\\&&&&\\ S &\longrightarrow & S^\sharp &\mbox{with}& S^\sharp(\xi):=Q\,S(Q^T\xi)\,Q^T\,.\end{array}\right\}\end{equation}
It can be shown (see \cite{MUNCHNEFF2016}) that 
\begin{equation}\label{curl-inc-inv-iso}\Curl_\xi\bigl[p^\sharp(\xi)\bigr]=Q\bigl(\Curl_xp(x)\bigr)Q^T\qquad\mbox{ and }\qquad \inc_\xi\bigl[\sym p^\sharp(\xi)\bigr]=Q\,\inc_x\bigl[\sym p(x)\bigr]\,Q^T\,.\end{equation}
Therefore, both our incompatibility measures are properly isotropic and therefore all our presented models, based on $\Curl p$ or $\inc(\sym p)$, respectively, are fully isotropic. 
\vskip.2truecm\noindent
A summary of the invariance conditions for infinitesimal gradient plasticity is presented in Table \ref{table:inv-inf-gradplast}.

 \begin{table}[h!]\footnotesize
\begin{center}\begin{tabular}{|ll|}\hline &\qquad\\
{\em Objectivity/Linearized frame-indifference:} &\hskip-.2truecm  $\left\{\begin{array}{lcl} \nabla u &\longrightarrow&\overline{A}+\nabla u,\\
p &\longrightarrow& \overline{A}+p\end{array}\right.\quad\overline{A}\in\so(3)$\\&\\
 {\em Linearized gauge-invariance:} & \hskip-.2truecm  $\left\{\begin{array}{lcl} \nabla u &\longrightarrow&\nabla u +\nabla\vartheta,\\
p &\longrightarrow& p+\nabla\vartheta\end{array}\right.\,\,\vartheta\in C^\infty(\mathbb{R}^3,\mathbb{R}^3)$\\&\\
 {\em Linearized micro-randomness:} & \hskip-.2truecm $\left.\begin{array}{lcl} \quad\,\, p\, &\longrightarrow& p +A(x)\end{array}\right.\quad A(x)\in\so(3)$ 
\\&\\
{\em Isotropy:}&\hskip-.2truecm $\left\{\begin{array}{lcll} h:\mathbb{R}^3\to\mathbb{R} &\longrightarrow & h^\sharp &\mbox{ with } h^\sharp(\xi):=h(Q^T\xi)\,,\\
\phi:\mathbb{R}^3\to\mathbb{R}^3&\longrightarrow & \phi^\sharp & \mbox{ with }\phi^\sharp(\xi):=Q\,\phi(Q^T\xi)\,,\\ S:\mathbb{R}^3\to\mathbb{R}^{3\times 3} &\longrightarrow & S^\sharp & \mbox{ with }S^\sharp(\xi):=Q\,S(Q^T\xi)\,Q^T\,\end{array}\right.\,Q\in\Oo(3)$ 
\\&\\
 \hline
\end{tabular}\end{center}\caption{\footnotesize Invariance conditions for infinitesimal gradient plasticity}\label{table:inv-inf-gradplast}\end{table}
\section{Some models of gradient plasticity with Kr\"{o}ner's
incompatibility tensor}\label{some models} Before we introduce and  analyze our "ideal"
model designed from the set of requirements presented in the
introduction, we found that it is more interesting to first present
those few models we first considered with an emphasis on the
difficulties and shortcomings of those models both from the
mechanical and mathematical points of view. Let us first make it
clear that the approach used to analyze those models as well as our
model in Section \ref{iso-hard-kroner} is trough a convex analytical
framework and variational inequalities developed in Han-Reddy
\cite{Han-ReddyBook} for classical plasticity and quite often used
for models of gradient plasticity (see \cite{DEMR1, REM, NCA,
EBONEFF, ENR2015}), as well.
\subsection{An irrotational model with linear kinematic hardening}\label{linkin-kroner} In this section, we present a model with linear kinematic
 hardening and Kr\"oner's incompatibility tensor where the plastic variable is symmetric i.e., a model with no plastic spin.
 The goal is to find the displacement field
$u$ and the infinitesimal plastic strain $\bvarepsilon_p$ in some
suitable function spaces
 such that the content of Table \ref{table:kinhardkroner1}
holds.
 \begin{table}[!h]\scriptsize
\begin{center}\begin{tabular}{|ll|}\hline &\qquad\\
 {\em Additive split of distortion:}& $\nabla u =e +p$\\
 {\em Additive split of strain:} &$\sym\nabla u=\bvarepsilon_e+\bvarepsilon_p$,\quad $\bvarepsilon_e:=\sym\,e\,\mbox{ and } \,\bvarepsilon_p:=\sym
 p$\\
{\em Equilibrium:} & $\mbox{Div}\,\sigma +f=0$ with
$\sigma=\C_{\mbox{\scriptsize{iso}}}\bvarepsilon_e$\\ {\em
Free energy:} &
$\frac12\,\langle\C_{\mbox{\scriptsize{iso}}}\bvarepsilon_e,\bvarepsilon_e\rangle+\frac12\,\mu\,
k_1\,\norm{\bvarepsilon_p}^2
+\frac12\,\mu\,\widehat{L}^4_c\norm{\inc(\bvarepsilon_p)}^2$\\&\\ {\em
Yield
condition:} & $\phi(\Sigma_E):=\norm{\dev\Sigma_E}-\yieldzero\leq0$\\
 {\em where } & $\Sigma_E:=\sigma+\Sigma_{\mbox{\scriptsize kin}} +\Sigma_{\mbox{\scriptsize inc}}$\,\,with
\,\,$\Sigma_{\mbox{\scriptsize
kin}}:=-\mu\, k_1 \bvarepsilon_p$,\\&
$\Sigma_{\mbox{\scriptsize
inc}}:=-\mu\,\widehat{L}^4_c\inc(\inc\bvarepsilon_p)\in\mbox{Sym}(3)$\\&
 \\{\em Dissipation inequality:} &
 $\dsize\int_\Omega\langle\Sigma_E,\dot{\bvarepsilon}_p\rangle\,
 dx\geq0$\\
 {\em Dissipation function:} &$\mathcal{D}(q):=\yieldzero \norm{q}$\\
 {\em Flow law in primal form:} &
 $\Sigma_E\in\partial \mathcal{D}(\dot{\bvarepsilon}_p)$\\
{\em Flow law in dual form:}
&$\dot{\bvarepsilon}_p=\lambda\,\dsize\frac{\dev\Sigma_E}{\norm{\dev\Sigma_E}},\quad\qquad
\lambda=\norm{\dot{\bvarepsilon}_p}$\\ {\em KKT conditions:}
&$\lambda\geq0$, \quad $\phi(\Sigma_E)\leq0$,
\quad $\lambda\,\phi(\Sigma_E)=0$\\
 {\em Boundary conditions for $\bvarepsilon_p$:} & to be specified\\
 {\em Function space for $\bvarepsilon_p$:} & $\bvarepsilon_p(t,\cdot)\in
 L^2(\Omega,\mbox{Sym}(3))$\,,\quad$\inc(\bvarepsilon_p)\in
 L^2(\Omega,\mbox{Sym}(3))$\\&\\
 \hline
\end{tabular}\end{center}\caption{\footnotesize The model with linear kinematic Prager-type hardening and Kr\"{o}ner's incompatibility tensor.  The boundary conditions
 on $\bvarepsilon_p$ cannot be specified from the structure of the model. The problem with  this model is the presence of the norm $\norm{\varepsilon_p}^2$ which is not {\bf linearized gauge-invariant}.  A careful analysis of the model in the spirit of Section
 \ref{iso-hard-kroner} will lead to single out the necessity of an $L^2$-control of $\Curl\bvarepsilon_p$, which is also not a priori controlled in the model.  Altogether,
well-posedness of the
 model  remains unclear. A modified version of this model
is presented in Table
\ref{table:kinhardkroner2}. This model is micro-random i.e., invariant w.r.t. $p\to p+A(x)$, $A(x)\in\so(3)$, but not linearized gauge-invariant, i.e., not invariant w.r.t.  $p\to p +\nabla\vartheta$, \,$\vartheta\in C^\infty(\mathbb{R}^3,\mathbb{R})$.}\label{table:kinhardkroner1}\end{table}
$\mbox{  }$\vskip-.5truecm
 \begin{table}[!tb]\scriptsize
\begin{center}\begin{tabular}{|ll|}\hline &\qquad\\
 {\em Additive split of distortion:}& $\nabla u =e +p$\\
 {\em Additive split of strain:} &$\sym\nabla u=\bvarepsilon_e+\bvarepsilon_p$,\quad $\bvarepsilon_e:=\sym\,e\,\mbox{ and } \,\bvarepsilon_p:=\sym
 p$\\
{\em Equilibrium:} & $\mbox{Div}\,\sigma +f=0$ with
$\sigma=\C_{\mbox{\scriptsize{iso}}}\bvarepsilon_e$\\ {\em
Free energy:} &
$\frac12\,\langle\C_{\mbox{\scriptsize{iso}}}\bvarepsilon_e,\bvarepsilon_e\rangle+\frac12\,\mu\,
k_1\,\norm{\bvarepsilon_p}^2 +\frac12\,\mu\,L_c^2\,\norm{\Curl\bvarepsilon_p}^2
+\frac12\,\mu\,\widehat{L}^4_c\norm{\inc(\bvarepsilon_p)}^2$\\&\\ {\em
Yield
condition:} & $\phi(\Sigma_E):=\norm{\dev\Sigma_E}-\yieldzero\leq0$\\
 {\em where } & $\begin{array}{lcl}\Sigma_E&:=&\sigma+\Sigma_{\mbox{\scriptsize kin}} +\Sigma_{\mbox{\scriptsize curl}} +\Sigma_{\mbox{\scriptsize inc}}\,,\qquad\Sigma_{\mbox{\scriptsize
kin}}:=-\mu\, k_1 \bvarepsilon_p\\
\Sigma_{\mbox{\scriptsize
curl}}&:=&-\mu\,L^2_c\sym\Curl\sym\Curl\bvarepsilon_p\in\mbox{Sym}(3)\\
\Sigma_{\mbox{\scriptsize
inc}}&:=&-\mu\,\widehat{L}^4_c\inc(\inc\bvarepsilon_p)\in\mbox{Sym}(3)\end{array}$
 \\&\\{\em Dissipation inequality:} &
 $\dsize\int_\Omega\langle\Sigma_E,\dot{\bvarepsilon}_p\rangle\,
 dx\geq0$\\
 {\em Dissipation function:} &$\mathcal{D}(q):=\yieldzero \norm{q}$\\&\\
 {\em Flow law in primal form:} &
 $\Sigma_E\in\partial \mathcal{D}(\dot{\bvarepsilon}_p)$\\
{\em Flow law in dual form:}
&$\dot{\bvarepsilon}_p=\lambda\,\dsize\frac{\dev\Sigma_E}{\norm{\dev\Sigma_E}},\quad\qquad
\lambda=\norm{\dot{\bvarepsilon}_p}$\\ {\em KKT conditions:}
&$\lambda\geq0$, \quad $\phi(\Sigma_E)\leq0$,
\quad $\lambda\,\phi(\Sigma_E)=0$\\
 {\em Boundary conditions for $\bvarepsilon_p$:} & $\bvarepsilon_p\times n|_{\partial\Omega}=0\quad\mbox{ and }\quad (\Curl\bvarepsilon_p)^T\times n|_{\partial\Omega}=0$\\
 {\em Function space for $\bvarepsilon_p$:} & $\bvarepsilon_p(t,\cdot)\in
 L^2(\Omega,\mbox{Sym}(3))$\,,\quad$\inc(\bvarepsilon_p)\in
 L^2(\Omega,\mbox{Sym}(3))$\\&\\
 \hline
\end{tabular}\end{center}

\caption{\footnotesize The model with linear kinematic Prager-type hardening and Kr\"{o}ner's incompatibility tensor.  The simple regularization term $\norm{\Curl\bvarepsilon_p}^2$ allows to  justify the  boundary conditions
 on $\bvarepsilon_p$.  The problem with  this model is the presence of the norms $\norm{\varepsilon_p}^2$  and $\norm{\Curl\bvarepsilon_p}^2$ which are both not {\bf linearized gauge-invariant}.   Altogether,
well-posedness of this modified
 model is obtained as in Section \ref{iso-hard-kroner}, with the difference that uniqueness of the weak solution is obtained without any additional regularity assumption. This model is micro-random i.e., invariant w.r.t. $p\to p+A(x)$, $A(x)\in\so(3)$, but not linearized gauge-invariant, i.e., not invariant w.r.t.  $p\to p +\nabla\vartheta$, $\vartheta\in C^\infty(\mathbb{R}^3,\mathbb{R})$.}\label{table:kinhardkroner2}\end{table}

\subsection{A fully isotropic model with isotropic hardening and plastic spin}
The model is completely described in Table \ref{table:isohard-fullyiso-new}.
\begin{table}[!t]\scriptsize\begin{center}
\begin{tabular}{|ll|}\hline &\qquad\\
 {\em Additive split of distortion:}& $\nabla u =e +p$,\quad $\bvarepsilon_e=\mbox{sym}\,e$,\quad $\bvarepsilon_p=\sym p$\\
{\em Equilibrium:} & $\mbox{Div}\,\sigma +f=0$ with
$\sigma=\C_{\mbox{\scriptsize{iso}}}\bvarepsilon_e$\\ {\em
Free energy:} &
$\dsize\frac12\,\langle\C_{\mbox{\scriptsize{iso}}}\bvarepsilon_e,\bvarepsilon_e\rangle+\frac12\,\mu\,
k_2\,|\widetilde{\gamma}_p|^2 +\,\,\frac12\,\mu\,L^2_c\norm{\Curl p}^2 $\\&\\
{\em Yield condition:} &$\phi(\Sigma_E,g):=\norm{\dev\Sigma_E} -(\yieldzero-g)\leq0$\quad where \quad$g:=-\mu\,k_2\widetilde{\gamma}_p$\\
 &$\begin{array}{lcl} \Sigma_E&:=&\sigma+\Sigma_{\mbox{\scriptsize
curl}}\,,\\
  \Sigma_{\mbox{\scriptsize curl}}&:=&-\,\mu\,L^2_c\,\Curl \Curl p\,\end{array}$
 \\&\\{\em Dissipation inequality:} &
 $\dsize\int_\Omega\bigl[\la\Sigma_E,\dot{p}\ra+g\,\dot{\widetilde{\gamma}}_p\ra\bigr]\,dx\geq0$\\
 {\em Dissipation function:} &$\mathcal{D}(\dot{\bvarepsilon}_p,\dot{\widetilde{\gamma}}_p):=\left\{\begin{array}{ll}\norm{\dot{p}} &\mbox{ if }
 \norm{\dot{p}}\leq\dot{\widetilde{\gamma}}_p,\\
 \infty &\mbox{ otherwise}\end{array}\right.$\\&\\
 {\em Flow law in primal form:} &
 $(\Sigma_E,g)\in\partial \mathcal{D}(\dot{\bvarepsilon}_p,\dot{\widetilde{\gamma}}_p)$\\
{\em Flow law in dual form:}
&$\dot{p}=\lambda\,\dsize\frac{\dev\Sigma_E}{\norm{\dev\Sigma_E}},
\quad\quad \dot{\widetilde{\gamma}}_p=\lambda=\norm{\dot{p}}$\\
{\em KKT conditions:} &$\lambda\geq0$, \quad $\phi(\Sigma_E,g)\leq0$,
\quad $\lambda\,\phi(\Sigma_E,g)=0$\\
 {\em Boundary conditions for $p$:} & $p\times n|_{\Gamma}=0$,\quad $\Curl p\times n|_{\partial\Omega\setminus\Gamma}=0$\\&\\
 \hline
\end{tabular}\caption{\footnotesize A fully isotropic model with isotropic hardening and plastic spin. The  model is fully isotropic since it is form-invariant under the $\sharp$-transformation defined  in (\ref{sharp-transfo}). The model is also linearized gauge-invariant i.e., invariant w.r.t.  $p\to p +\nabla\vartheta$, $\vartheta\in C^\infty(\mathbb{R}^3,\mathbb{R})$, but not micro-random i.e., not invariant w.r.t. $p\to p+A(x)$,  $A(x)\in\so(3)$.}\label{table:isohard-fullyiso-new}\end{center}\end{table}

\subsection{An irrotational model with isotropic
hardening}\label{iso-hard-kroner} In this section we discuss a variant
of the previous model with linear kinematic hardening replaced
by isotropic hardening.  The new
model will be invariant under Linear Referential Isotropy ({\bf LRIso}),
Linear Micro-Random ({\bf LMR}), Linear Gauge-Invariance ({\bf LGI}), Linear Elastic
Objectivity, Linear Elastic Isotropy. \subsubsection{Derivation of the model}
{\bf The balance equation.} The
conventional macroscopic force balance leads to the equation of
equilibrium
\begin{equation}
\Div \sigma + f = 0 \label{equil}
\end{equation}
in which $\sigma$ is the infinitesimal symmetric Cauchy stress and
$f$ is the body force.\vskip.1truecm\noindent {\bf Constitutive
equations.} The constitutive equations are obtained from a free
energy imbalance together with a flow law that characterizes plastic
behaviour. The total strain $\bvarepsilon$ is additively decomposed
into elastic and plastic components $\bvarepsilon_e$ and
$\bvarepsilon_p$, so that
\begin{equation}
\bvarepsilon = \bvarepsilon_e + \bvarepsilon_p \label{strain}
\end{equation}
with the plastic strain incapable of sustaining volumetric changes;
that is, $\mbox{tr} \,\bvarepsilon_p = 0\,.$\\
The strain-displacement relation is given by
\begin{equation}\label{strain-displ}
\bvarepsilon=\sym\nabla u=\frac12(\nabla u+\nabla
u^T)\,.\end{equation}
 {\bf Free energy density:} 
 In this model the free-energy density is considered in the  additively separated  form
\begin{eqnarray}\label{free-eng-iso}
\mbox{W}(u,\varepsilon_p, \inc\varepsilon_p,\gamma_p):
&=&\underbrace{\mbox{W}_e(\bvarepsilon_e)}_{\mbox{\footnotesize elastic energy}}\,\,+\,\,
\underbrace{\mbox{W}_{\mbox{\scriptsize
inc}}(\inc\varepsilon_p)}_{\begin{array}{c}\mbox{\footnotesize mesoscopic}\\
\mbox{\footnotesize incompatibility}\end{array}}\,\, +\,\,\underbrace{\mbox{W}_{\mbox{\scriptsize
 iso}}(\gamma_p)}_{\begin{array}{c}\mbox{\footnotesize isotropic hardening}\\
 \mbox{\footnotesize energy (SSD-related)}\end{array}}\,,
 \end{eqnarray} where
 \begin{equation}\label{eng-polyg1}\left.\begin{array}{lll} \mbox{W}_{\mbox{\footnotesize e}}(\bvarepsilon_e)&:=&\dsize\frac12\,\la
\bvarepsilon_e,\C_{\mbox{\scriptsize{iso}}}\bvarepsilon_e\ra =\mu\,\norm{\sym\nabla u-\varepsilon_p}^2+\frac12\,\lambda\tr[\nabla u -\varepsilon_p]^2\,, \\\\
 \mbox{W}_{\mbox{\scriptsize
inc}}(\inc\varepsilon_p)&:=&\dsize\frac12\,\mu\,\widehat{L}^4_c\norm{\inc\varepsilon_p}^2=\frac12\mu\,\widehat{L}^4_c\norm{\Curl[(\Curl\varepsilon_p)^T]  }^2,\\\\
 \mbox{W}_{\mbox{\scriptsize
 iso}}(\gamma_p)&:=&\dsize\frac12\,\mu\,k_2\,|\gamma_p|^2\,,\end{array}\right\}\end{equation}
and $\lambda,\,\mu$ are the Lam\'e moduli  with $\mu>0$ and $3\lambda+2\mu>0$, $\widehat{L}_c\geq0$ is an energetic length scale and $k_2\geq0$ is a  positive non-dimensional isotropic hardening constant,
$\gamma_p$ is the isotropic hardening variable (the accumulated equivalent plastic strain).\vskip.3truecm
From the local free energy inbalance 
\begin{eqnarray*}
\frac{\rm d}{\rm dt} \mbox{W}\leq\la\sigma,\nabla u_t\ra =\la\sigma,\dot{\bvarepsilon}\ra&\Leftrightarrow &\frac{\rm d}{\rm dt}\mbox{W} - \la\sigma,\dot{\varepsilon}_e\ra - \la\sigma,\dot{\varepsilon}_p\ra\leq0\\
&\Leftrightarrow &\frac{\rm d}{\rm dt}{\rm W}(u,\varepsilon_p,\mbox{D}^2\varepsilon_p,\gamma_p)\leq0\quad\mbox{ for $u$ fixed}\end{eqnarray*}
where the second equivalence is obtained using arguments from thermodynamics which give the elasticity
relation
\begin{equation}
\sigma = {\C}_{\mbox{\scriptsize{iso}}}\bvarepsilon_e=2\mu\,
(\sym\nabla u-\bvarepsilon_p)+\lambda\, \tr(\sym\nabla
u-\bvarepsilon_p)\id\,. \label{elasticlaw}
\end{equation}
Therefore, we get
\begin{equation}\label{diss-ineq-loc1}
\la\sigma,\dot{\varepsilon}_p\ra 
-\mu\,\widehat{L}^4_c\la\inc\varepsilon_p,\inc\dot{\varepsilon}_p\ra
-\mu\,k_2\,\gamma_p\,\dot{\gamma}_p\geq0\,.\end{equation}

Now, integrating (\ref{diss-ineq-loc1}),  we arrive at

\begin{eqnarray}\label{constr-dissip}\nonumber0&\leq& \int_{\Omega}\Bigl[\la\sigma,\dot{\bvarepsilon}_p\ra -\mu\,
\widehat{L}_c^4\,\la\,\underbrace{\inc(\inc\bvarepsilon_p)}_{\in\mbox{Sym}(3)},\dot{\bvarepsilon}_p\ra -\mu\,k_2\,\gamma_p\,\dot{\gamma}_p-\sum_{i=1}^3\div\Bigl(\mu\,
\widehat{L}_c^4\,\dot{\bvarepsilon}_{p\,i}\times\bigl[\Curl\inc\bvarepsilon_p\bigr]^T_i\Bigr)\\
&&\nonumber \hskip1.5truecm
 -\sum_{i=1}^3\div\Bigl(\mu\,
\widehat{L}_c^4\,\bigl[\Curl\dot{\bvarepsilon}_p\bigr]^T_i\times\bigl[\inc\bvarepsilon_p\bigr]_i\Bigr)\Bigr]\,dx\\
&=& \int_{\Omega}\Bigl[\la\sigma-\mu\, \widehat{L}_c^4\,\inc(\inc
\bvarepsilon_p),\dot{\varepsilon}_p\ra-\mu\,k_2\,\gamma_p\,\dot{\gamma}_p\Bigr]dx\\
&&\nonumber\hskip1.truecm-\sum_{i=1}^3\mu\,\widehat{L}_c^4\int_{\partial\Omega}\la\dot{\bvarepsilon}_{p\,i}\times\bigl[\Curl\inc\bvarepsilon_p\bigr]^T_i,n\ra\,dS-\sum_{i=1}^3\mu\,
\widehat{L}_c^4\int_{\partial\Omega}\la\bigl[\Curl\dot{\bvarepsilon}_p\bigr]^T_i\times\bigl[\inc\bvarepsilon_p\bigr]_i,n\ra\,dS\,.\end{eqnarray}

In order to obtain a  global reduced dissipation
inequality one needs to choose suitable boundary conditions for which the
two equations below are satisfied
\begin{eqnarray}
\label{lin-insul1}
\sum_{i=1}^3\int_{\partial\Omega}\la\dot{\bvarepsilon}_{p\,i}\times n,\bigl[\Curl\inc\bvarepsilon_p\bigr]^T_i\ra\,dS&=&0\,.\\&&\nonumber\\
\label{lin-insul2}\sum_{i=1}^3\int_{\partial\Omega}\la\bigl[\Curl\dot{\bvarepsilon}_p\bigr]^T_i\times n,
\bigl[\inc\bvarepsilon_p\bigr]_i\ra\,dS &=&0\,.\end{eqnarray} The
simplest lower order boundary conditions to satisfy (\ref{lin-insul1}) and (\ref{lin-insul2}) are
\begin{equation}\label{boundary-conditions}
 \varepsilon_p\times n|_{\partial\Omega}=0\qquad\mbox{ and }\qquad
\bigl[\Curl\varepsilon_p\bigr]^T\times n|_{\partial\Omega}=0\,.
\end{equation}
Other possible boundary conditions to satisfy the equations (\ref{lin-insul1}) and (\ref{lin-insul2})    are given in the table below. 
\begin{table}[h!]\small\begin{center}
\begin{tabular}{|l|l|}\hline &\\ 
 {\bf  boundary conditions for (\ref{lin-insul1})} &  {\bf boundary conditions for (\ref{lin-insul2})}\\
 \hline&\\
$ [\Curl\inc\bvarepsilon_p\bigr]^T\times n|_{\partial\Omega}=0$ & $(\inc\varepsilon_p)\times n|_{\partial\Omega}=0$\\&\\
 \hline&\\
 $\varepsilon_p\times n|_{\Gamma}=0\quad\mbox{ and }\quad
\bigl[\Curl\inc\varepsilon_p\bigr]^T\times n|_{\partial\Omega\setminus\Gamma}=0$ &
$\bigl[\Curl\varepsilon_p\bigr]^T\times n|_{\Gamma}=0\quad\mbox{ and }\quad \inc\bvarepsilon_p\times n|_{\partial\Omega\setminus\Gamma}=0$\,.
\\&\\
\hline
 
\end{tabular}\caption{\footnotesize Possible boundary conditions for (\ref{lin-insul1}) and (\ref{lin-insul2}) to be satisfied. }\label{table:bc}\end{center}\end{table}

However, 
these boundary conditions cannot be mathematically  justified from the free-energy density $\mbox{W}$ considered so far: both terms in (\ref{lin-insul1}) and (\ref{lin-insul2}) are not automatically well-defined as boundary traces. In fact, one needs to show that $\bvarepsilon_p\in \mbox{H}(\mbox{Curl})$ and $(\Curl\bvarepsilon_p)^T\in \mbox{H}(\mbox{Curl})$. This information is missing from the energy. We only know that $\bvarepsilon_p\in\mbox{L}^2$ (due to isotropic hardening) and $\inc\bvarepsilon_p=\Curl[(\Curl\bvarepsilon_p)^T]\in\mbox{L}^2$. The missing piece of information to proceed  is $\Curl\bvarepsilon_p\in\mbox{L}^2$.

So, one needs to modify the model by adding a new regularizing term in the free-energy density $W$, which is physically meaningful in the sense that it does satisfy some invariance properties. The unmodified model is summarized  in Table \ref{table:isohard-inc}.

\begin{table}[h!]\footnotesize\begin{center}
\begin{tabular}{|ll|}\hline &\qquad\\
 {\em Additive split of strain:} &$\sym\nabla u=\bvarepsilon_e+\bvarepsilon_p$\\
{\em Equilibrium:} & $\mbox{Div}\,\sigma +f=0$ with
$\sigma=\C_{\mbox{\scriptsize{iso}}}\bvarepsilon_e$\\&\\ {\em
Free energy:} &
$\dsize\frac12\,\langle\C_{\mbox{\scriptsize{iso}}}\bvarepsilon_e,\bvarepsilon_e\rangle+\frac12\,\mu\,
k_2\,|\gamma_p|^2  +\,\frac12\,\mu\,\widehat{L}^4_c\norm{\inc\bvarepsilon_p}^2$\\&\\
{\em Yield condition:} &$\phi(\Sigma_E,g):=\norm{\dev\Sigma_E} -(\yieldzero-g)\leq0$\quad where \quad$g:=-\mu\,k_2\gamma_p$\\
 &$\begin{array}{lcl} \Sigma_E&:=&\sigma+\Sigma_{\mbox{\scriptsize
 inc}}\,,\\
\Sigma_{\mbox{\scriptsize
inc}}&:=&-\mu\,\widehat{L}^4_c\,\inc\inc\bvarepsilon_p\,=\,\Curl\Bigl(\bigl[\Curl\Curl[(\Curl \bvarepsilon_p)^T] \bigr]^T\Bigr)\end{array}$
\\& \\{\em Dissipation inequality:} &
 $\dsize\int_\Omega\bigl[\la\Sigma_E,\dot{\bvarepsilon}_p\ra+g\,\dot{\gamma}_p\ra\bigr]\,dx\geq0$\\
 {\em Dissipation function:} &$\mathcal{D}(\dot{\bvarepsilon}_p,\dot{\gamma}_p):=\left\{\begin{array}{ll}\yieldzero\, \norm{\dot{\bvarepsilon}_p} &\mbox{ if }
 \norm{\dot{\bvarepsilon}_p}\leq\dot{\gamma}_p,\\
 \infty &\mbox{ otherwise}\end{array}\right.$\\&\\
 {\em Flow law in primal form:} &
 $(\Sigma_E,g)\in\partial \mathcal{D}(\dot{\bvarepsilon}_p,\dot{\gamma}_p)$\\
{\em Flow law in dual form:}
&$\dot{\bvarepsilon}_p=\lambda\,\dsize\frac{\dev\Sigma_E}{\norm{\dev\Sigma_E}},
\quad\quad \dot{\gamma}_p=\lambda=\norm{\dot{\bvarepsilon}_p}$\\
{\em KKT conditions:} &$\lambda\geq0$, \quad $\phi(\Sigma_E,g)\leq0$,
\quad $\lambda\,\phi(\Sigma_E,g)=0$\\
 {\em Boundary conditions for $\bvarepsilon_p$:} & $\bvarepsilon_p\times n|_{\partial\Omega}=0\,,\quad (\Curl \bvarepsilon_p)^T\times n|_{\partial\Omega}=0$\\&\\
 \hline
\end{tabular}\caption{\footnotesize The irrotational model with isotropic hardening and Kr\"{o}ner's incompatibility tensor.  The terms in the free-energy density are not enough to guarantee the minimum regularity required, i.e., $\bvarepsilon_p\in \mbox{H}(\mbox{Curl})$ and $(\Curl\bvarepsilon_p)^T\in \mbox{H}(\mbox{Curl})$,  in order to justify mathematically the boundary conditions $\bvarepsilon_p\times n|_{\partial\Omega}=0$ and $ (\Curl \bvarepsilon_p)^T\times n|_{\partial\Omega}=0$. }\label{table:isohard-inc}\end{center}\end{table}

\vskip.2truecm\noindent
We will consider the additional term 
\begin{equation}\label{new-term1}
\mbox{W}_{\mbox{\scriptsize
 curl}}(\Curl p):=\dsize\frac12\,\mu\, L_c^2\,\norm{\dev\sym\Curl \bvarepsilon_p}^2\,,
\end{equation}
which is motivated in the following section.\vskip.2truecm\noindent
\subsubsection{Conformal gauge-invariance - the regularization term $\mathbf\dev\sym\Curl\bvarepsilon_p$}
 We will see subsequently  that the model with the regularizing term $\dev\sym\Curl\bvarepsilon_p$  allows for a mathematical existence proof. However, what about the invariance conditions, notably gauge-invariance?\\
 It is easy to see that 
$$ p\to \dev\sym\Curl\bvarepsilon_p=\dev\sym\Curl\sym p$$ is micro-random while it is not linear gauge-invariant, i.e.,
\begin{equation}\label{invariance-new1}\begin{array}{rcl} p &\to& p+\nabla\vartheta\\
\dev\sym\Curl\sym(\nabla\vartheta +p) &\neq &\dev\sym\Curl\sym p,\quad\forall\vartheta\in C^1(\Omega,\mathbb{R}^3)\,.\end{array}\end{equation}
Let us now determine those mappings $\vartheta:\mathbb{R}^3\to\mathbb{R}^3$  which are still ``allowed'' for gauge-invariance, in the sense that 
$$\dev\sym\Curl\sym(\nabla\vartheta +p) \,=\,\dev\sym\Curl\sym p\,.$$
 Automatically, these mappings satisfy the identity $\inc(\sym\nabla\vartheta)=0$. 
Moreover, by linearity we should have 
\begin{equation}\label{inv-conf1}\dev\sym\Curl\sym\nabla\vartheta=0\,.\end{equation} Since, however, $\tr(\Curl S)=0$ for all smooth symmetric tensor fields $S\in\mbox{Sym}(3)$ (see (\ref{trace-curl-sym2}) in the  appendix), the latter is equivalent to
$$\sym\Curl\sym\nabla\vartheta=0\,.$$
This implies that for  some non-constant skew-symmetric tensor field $A:\Omega\to\so(3)$ we have 
\begin{equation}\label{invariance-new2}\Curl\sym\nabla\vartheta(x)=A(x)\quad\Leftrightarrow\quad (\Curl\sym\nabla\vartheta)^T=-A(x)\,.\end{equation} 
Taking the Curl on both sides leads to 
\begin{equation}\label{inv-new1}\Curl[(\Curl \sym\nabla \vartheta)^T]=-\Curl A(x)\quad\Leftrightarrow\quad 0\,=\,\inc(\sym\nabla\vartheta)\,=-\,\Curl A(x)\,.\end{equation}

Thus, $A(x)=\overline{A}$ is a constant skew-symmetric matrx,  according to an observation in  \cite{MUNCHNEFF2016}. Reinserting into (\ref{invariance-new2}), we must have \begin{equation}\label{reg-add-eq1}\Curl\sym\nabla\vartheta(x)=\overline{A}\,.\end{equation}
We observe that (see \cite{NGMPR2014})
$$\Curl(\zeta(x_1,x_2,x_3)\cdot\id)=\left(\begin{array}{ccc} 0 & -\zeta_{,3} & \zeta_{,2}\\
 \zeta_{,3}&0 & -\zeta_{,1}\\
 -\zeta_{,2}&\zeta_{,1}&0\end{array}\right)\, \quad \zeta:\mathbb{R}^3\to\mathbb{R}\,$$ 
and with $\zeta(x_1,x_2,x_3)=a\,x_1+b\,x_2+c\,x_3 $, we obtain 
$$\Curl((ax_1+bx_2+cx_3)\cdot\id)=\left(\begin{array}{ccc} 0 & -c & b\\
 c&0 & -a\\
 -b&a&0\end{array}\right)\,.$$ 
 
 Hence, a solution to (\ref{reg-add-eq1}) can be obtained in the format 
 \begin{equation}\label{invariance-new3}\sym\nabla\vartheta (x_1,x_2,x_3)=\zeta(x_1,x_2,x_3)\cdot\id\,.\end{equation} On taking again the deviatoric part of the latter  we arrive at 
 \begin{equation}\label{reg-add-eq2}
 \dev\sym\nabla\vartheta(x_1,x_2,x_3)=0\,.\end{equation}
 This is equivalent to \begin{equation}\label{inv-new3}\nabla\vartheta (x)=\zeta(x_1,x_2,x_3)\cdot\id + A(x_1,x_2,x_3)\quad\mbox{ with }\quad A:\mathbb{R}^3\to\so(3)\,.\end{equation}
 The solution to (\ref{reg-add-eq2}) can be given in closed form. In fact, taking Curl on both sides of (\ref{inv-new3}), together with the fact that $\Curl(\zeta(x_1,x_2,x_3)\cdot\id)\in\so(3)$ one gets $A(x_1,x_2,x_3)=\widehat{A}\in\so(3)$ constant skew-symmetric matrix. Also, using the operators {\bf axl} and {\bf anti} defined in (\ref{axl-anti1}), a general solution to (\ref{reg-add-eq2}) is obtained in the form
 \begin{equation}\label{inv-closed-new1}\phi_c(x):=\frac12\,\Bigl(2\la\axl\bigl(\widehat{W}\bigr),x\ra\,x-\axl\bigl(\widehat{W}\bigr)\,\norm{x}^2\Bigr)+\bigl[\widehat{\zeta}\cdot\id +\widehat{A}\bigr]\,x+\widehat{\eta}\,,\end{equation}
 where $ \widehat{A},\,\widehat{W}\in\so(3)$ are arbitrary constant skew-symmetric matrices and $\widehat{\zeta},\,\widehat{\eta}\in\mathbb{R}^3$ are arbitrary constant vectors.\\
 The mappings in (\ref{inv-closed-new1}) are called {\bf infinitesimal conformal mappings} $\phi_c$ (see \cite{NJ2009}).
 The mappings $x\to\phi_c(x)$ locally preserve the shape of infinitesimal cubes but are globally inhomogeneous.
 
\begin{figure}[h!]
\centering
\vskip-8truecm
\includegraphics[width=1.\textwidth, height=1.\textheight]{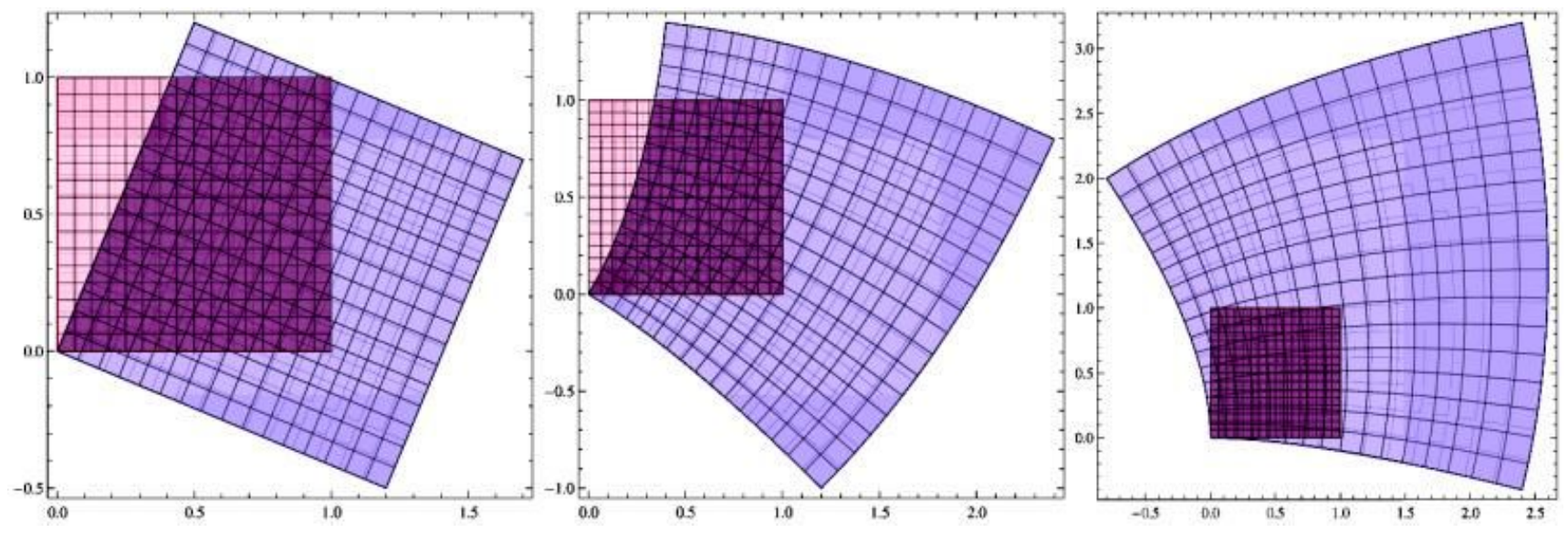}
\vskip-8.5truecm\caption{\footnotesize Infinitesimal conformal mappings $\phi_c:\mathbb{R}^3\to\mathbb{R}^3$ which locally leave shapes invariant: a prototype elastic deformation in the sense that the corresponding stress deviator $\dev\sigma(\nabla\phi_c)=0$. Shown is the coarse grid
deformation. The picture is extracted from \cite{NJ2009}.}\label{fig:inf-conf-mapping}
\end{figure}
  
  If we consider $x\to\phi_c(x)$ as elastic displacement, then, according to the von Mises $J_2$-criterion, these mappings alone  never lead to plasticity since 
  $$\dev\sigma=\dev\bigl(2\,\mu\,\sym\nabla\phi_c+\lambda\,\tr(\nabla \phi_c\cdot\id\bigr)=2\,\mu\dev\sym\nabla\phi_c=0\,.$$
 
 Gathering our findings, we have obtained that the regularization term (\ref{new-term1}) is invariant w.r.t. the infinitesimal conformal group and  infinitesimal conformal mappings $\phi_c$ do not induce irreversible processes.
 
 There is still another solution to 
 \begin{equation}\label{reg-add-eq3}\sym\Curl\sym\nabla\vartheta=0\,.\end{equation}
  Clearly, (\ref{reg-add-eq3}) will be satisfied also if  already 
  $$\Curl\sym\nabla\vartheta=0\,,$$ which in turn is satsfied for $\sym\nabla\vartheta=\nabla v\in\mbox{Sym}(3)$, with $v:\mathbb{R}^3\to\mathbb{R}^3$. Such a vector can be taken as $v=\nabla h(x_1,x_2,x_3)$ with any scalar function $h:\mathbb{R}^3\to\mathbb{R}$. Then, (\ref{reg-add-eq3}) is satisfied. Thus, another solution to (\ref{reg-add-eq3})  is given by 
  $$\vartheta(x_1,x_2,x_3)=\nabla h(x_1,x_2,x_3),\quad h:\mathbb{R}^3\to\mathbb{R}\,.$$
  Altogether, solutions to (\ref{reg-add-eq3}) are represented by 
  $$\boxed{\vartheta(x)= \underbrace{\phi_c(x)}_{\mbox{\scriptsize``conformal''}}\,+\underbrace{\nabla h}_{\mbox{\scriptsize  ``potential''}}}$$ as the new invariance group. We collect our finding in the following theorem.
  \begin{theorem}\label{null-devsymcurl}{\rm \bf [Nullspace of dev\,sym\,Curl\,sym\,Grad]}\\
  The nullspace of the operator  $\dev\sym\Curl\sym \Grad$ is given by 
  $$\vartheta(x)=\frac12\,\Bigl(2\la\axl\bigl(\widehat{W}\bigr),x\ra\,x-\axl\bigl(\widehat{W}\bigr)\,\norm{x}^2\Bigr)+\bigl[\widehat{\zeta}\cdot\id +\widehat{A}\bigr]\,x+\widehat{\eta} +\nabla h(x)\,,$$
 where $ \widehat{A},\,\widehat{W}\in\so(3)$ are arbitrary constant skew-symmetric matrices, $\widehat{\zeta},\,\widehat{\eta}\in\mathbb{R}^3$ are arbitrary constant vectors and $h:\mathbb{R}^3\to\mathbb{R}$ is any scalar function.

  \end{theorem}
  
  It is remarkable, that the seemingly similar regularization term 
  $\norm{\Curl\sym p}^2=\norm{\Curl\bvarepsilon_p}^2$ only allows for invariance under ``potential'' mappings $\vartheta=\nabla h$. 
  \vskip.2truecm\noindent

\begin{figure}[h!]\tiny
\centering

\includegraphics[width=0.5\textwidth, height=0.5\textheight]{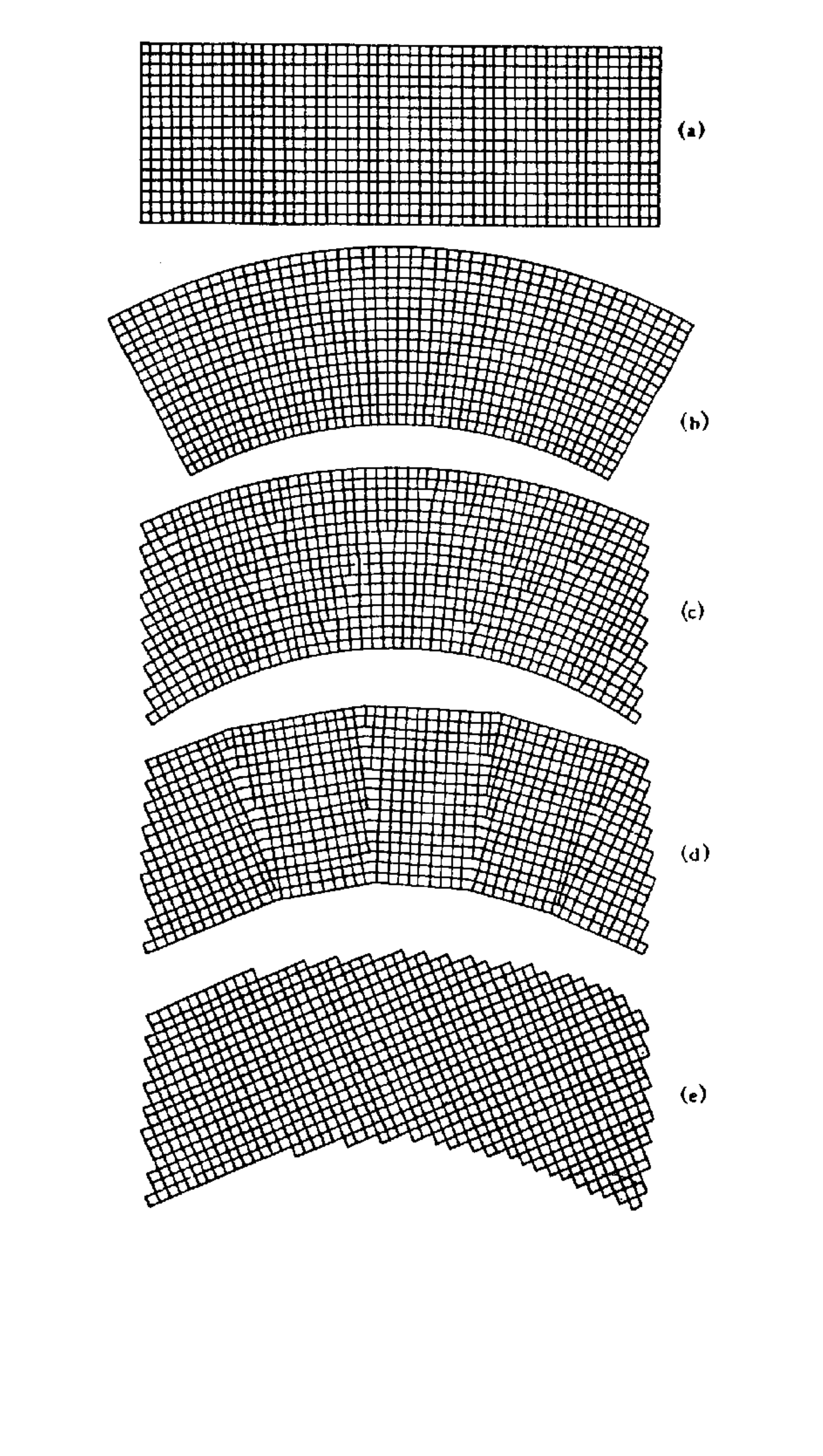}
\vskip-1.5truecm\caption{\footnotesize Five states of a single crystal. (a) Unstrained. (b) Elastically bent. (c) Plastically bent. (d) Polygonized. \quad (e) Recrystallized. The picture was extracted from Nabarro \cite{Nabarro1967}. 
 For the polygonized crystal, we observe large structures which are rotated against each other with a ``zone'' separating those blocks.\,\,Using plainly $\Curl p$ as the underlying defect measure would energetically penalize these configurations. Therefore it seems appropriate to use $\inc(\sym p)$ as a ``weaker'' defect measure which would allow for low energy configurations like that in (d).
}\label{fig:polygo}
\end{figure}

   In order to be able to describe polygonization (see Figure \ref{fig:polygo}(d)), the plasticity model should energetically favour configurations in which there are blocks of many homogeneous rotations.\\
   In this respect, the new term $\frac12\,\mu\,L_c^2\norm{\dev\sym\Curl\sym p}^2=\frac12\,\mu\,L_c^2\norm{\sym\Curl\sym p}^2$ energetically favours those configurations, which locally only rotate. The generated natural second order backstress will be of the type
   $$\Sigma_{\mbox{\scriptsize curl}}=\mu\,L_c^2\,\sym\Curl(\sym\Curl\bvarepsilon_p)\,.$$
   Now, looking at the  invariance of the energy for which $\sym\Curl(\sym\nabla\vartheta)=0$ versus the invariance of the backstress in the strong formulation for which $\sym\Curl(\sym\Curl(\sym\nabla\vartheta)=0$, it is clear that the invariance of the energy implies the  invariance of backstress, but not vice-versa.
   \begin{remark}\label{rem-on-new-term}{\rm 
   Note that the mapping $p\to\dev\sym\Curl\sym p=\sym\Curl\bvarepsilon_p$ does not have any geometric meaning connected to the incompatibility of the plastic distortion $p$ like $\Curl p$ or connected to the incompatibility of the plastic strain tensor $\bvarepsilon_p=\sym p$ like $\inc(\sym p)$. The simpler term $\Curl\sym p$ has been used by Gurtin and Anand \cite{GURTAN2005} as the only energetic contribution  in their irrotational gradient plasticity  model.}\end{remark}
\subsubsection{Derivation of the modified model}
Now with the additional term 
$\mbox{W}_{\mbox{\scriptsize
 curl}}(\Curl p):=\dsize\frac12\,\mu\, L_c^2\,\norm{\dev\sym\Curl \bvarepsilon_p}^2=\dsize\frac12\,\mu\, L_c^2\,\norm{\sym\Curl \bvarepsilon_p}^2$ in the free-energy density $\mbox{W}$, if we repeat the derivation above starting from the free-energy imbalance, we get
\begin{eqnarray}\label{diss-ineq-new1}
0&\leq& \int_{\Omega}\Bigl[\la\sigma-\mu\,L^2_c\,\sym\Curl(\sym\Curl\varepsilon_p)-\mu\, \widehat{L}_c^4\,\inc(\inc
\bvarepsilon_p),\dot{\varepsilon}_p\ra-\mu\,k_2\,\gamma_p\,\dot{\gamma}_p\Bigr]dx\\
&&\nonumber\hskip1.truecm-\sum_{i=1}^3\mu\,L_c^2\int_{\partial\Omega}\la\dot{\bvarepsilon}_{p\,i}\times n,(\sym\Curl\bvarepsilon_p)_i\ra\,dS -\sum_{i=1}^3\mu\,\widehat{L}_c^4\int_{\partial\Omega}\la\dot{\bvarepsilon}_{p\,i}\times n,\bigl[\Curl\inc\bvarepsilon_p\bigr]^T_i\ra\,dS\\
&&\nonumber\hskip1.truecm-\sum_{i=1}^3\mu\,
\widehat{L}_c^4\int_{\partial\Omega}\la\bigl[\Curl\dot{\bvarepsilon}_p\bigr]^T_i\times n,\bigl[\inc\bvarepsilon_p\bigr]_i\ra\,dS\\
&=&\int_{\Omega}\Bigl[\la\Sigma_E,\dot{\varepsilon}_p\ra +g\dot{\gamma}_p\Bigr]dx
-\sum_{i=1}^3\mu\,L_c^2\int_{\partial\Omega}\la\dot{\bvarepsilon}_{p\,i}\times n,(\sym\Curl\bvarepsilon_p)_i\ra\,dS\\
&&\nonumber\hskip1.truecm-\sum_{i=1}^3\mu\,\widehat{L}_c^4\int_{\partial\Omega}\la\dot{\bvarepsilon}_{p\,i}\times n,\bigl[\Curl\inc\bvarepsilon_p\bigr]^T_i\ra\,dS-\sum_{i=1}^3\mu\,
\widehat{L}_c^4\int_{\partial\Omega}\la\bigl[\Curl\dot{\bvarepsilon}_p\bigr]^T_i\times n,\bigl[\inc\bvarepsilon_p\bigr]_i\ra\,dS
\end{eqnarray}
where 
\begin{eqnarray*}
\Sigma_E &:=&\sigma +\Sigma_{\mbox{\scriptsize
 curl}}+\Sigma_{\mbox{\scriptsize
inc}}\,,\qquad g=-\,\mu\,k_2\,\gamma_p\,,\\
\Sigma_{\mbox{\scriptsize
 curl}}&:=&-\,\mu\,L^2_c\,\sym\Curl(\sym\Curl \varepsilon_p)\,,\quad\qquad\mbox{(second order nonlocal backstress)}\\
 \Sigma_{\mbox{\scriptsize
  inc}}&:=&-\,\mu\,\widehat{L}^4_c\,\inc(\inc
\bvarepsilon_p)\\
&=&-\mu\,\widehat{L}^4_c\,\Curl\Bigl(\Bigl[\Curl\Curl\bigl[\Curl\varepsilon_p\bigr]^T\Bigr]^T\Bigr)\quad\,\mbox{(fourth order nonlocal backstress)}\,.\end{eqnarray*}
Now assuming again the simplest lower order boundary conditions
\begin{equation}\label{bc-new}
 \varepsilon_p\times n|_{\partial\Omega}=0\quad\mbox{ and }\quad
\bigl[\Curl\varepsilon_p\bigr]^T\times n|_{\partial\Omega}=0
\end{equation}
which will be clearly defined as Sobolev traces through a choice of a suitable function space for the plastic strain variable $\bvarepsilon_p$, 
will guarantee the insulation type conditions 
\begin{eqnarray}
\sum_{i=1}^3\int_{\partial\Omega}\la\dot{\bvarepsilon}_{p\,i}\times n,(\sym\Curl
\bvarepsilon^p)_i\ra\,dS &=&0\,,\label{lin-insul0}\\&&\nonumber\\
\label{lin-insul1-new}
\sum_{i=1}^3\int_{\partial\Omega}\la\dot{\bvarepsilon}_{p\,i}\times n,\bigl[\Curl\inc\bvarepsilon_p\bigr]^T_i\ra\,dS&=&0\,,\\&&\nonumber\\
\label{lin-insul2-new}\sum_{i=1}^3\int_{\partial\Omega}\la\bigl[\Curl\dot{\bvarepsilon}_p\bigr]^T_i\times n,
\bigl[\inc\bvarepsilon_p\bigr]_i\ra\,da &=&0\,,\end{eqnarray}
from which we obtain the global reduced dissipation
inequality
\begin{equation}\label{global-reduced-dissip}\int_{\Omega}\bigl[\la\Sigma_E,\dot{\varepsilon}_p\ra +g\,\dot{\gamma}_p\bigr]\,dx\,\geq\,0\,.\end{equation}
\vskip.2truecm\noindent
{\bf The flow
law:} 
  We consider 
the set of admissible (elastic) generalized stresses 
\begin{equation}\label{admin-stresses}
\mathcal{E}:=\{(\Sigma_E,g)\,\,|\,\,
\norm{\dev\Sigma_E}-(\sigma_0-g)\leq0,\quad g\leq0\}\,,\end{equation} whose interior Int$(\mathcal{E})$ is the elastic domain while its boundary $\partial\mathcal{E}$ is the yield surface. The constant $\yieldzero$ is the initial yield stress of the material. The flow law in
its primal form reads as follows:
\begin{equation}\label{flow-primal-iso}
(\Sigma_E,g)\in\partial\mathcal{D}(\dot{\bvarepsilon}_p,\dot{\gamma}_p)\end{equation}  where
\begin{eqnarray}\label{diss-funt-iso}
\nonumber\mathcal{D}(q,\beta)&:=&\sup\{\la\Sigma_E,q\ra+g\beta\,\,\,|\,\,\,(\Sigma_E,g)\in\mathcal{E}\}\\
&=&\nonumber\sup\{\la\Sigma_E,q\ra+g\,\beta\,\,\,|\,\,\,\norm{\Sigma_E}\leq\sigma_0-g,\quad g\leq0\}\\
&=& \left\{\begin{array}{ll} \yieldzero\,\norm{q}
&\mbox{
if } \norm{q}\leq\beta,\\
\infty & \mbox{ otherwise.}\end{array}\right. \end{eqnarray} 
 Here, $\partial\mathcal{D}(\dot{\Gamma}_p)$ denotes the
 subdifferential of the function $\mathcal{D}$ at $\dot{\Gamma}_p$.
 That is,
 \begin{equation}\label{primal-flow-law}(\Sigma_E,g)\in\partial\mathcal{D}(\dot{\bvarepsilon}_p,\dot{\gamma}_p)\quad \Leftrightarrow\quad
 \la\Sigma_E,q-\dot{\bvarepsilon}_p\ra+g(\beta-\dot{\gamma}_p)\leq\mathcal{D}(q,\beta)-\mathcal{D}(\dot{\bvarepsilon}_p,\dot{\gamma}_p)\quad\forall\,(q,\beta)\,.\end{equation}
Now using convex analysis, we get
\begin{equation}\label{equivflowlaw}
(\Sigma_E,g)\in\partial\mathcal{D}(\dot{\bvarepsilon}_p,\dot{\gamma}_p)\quad
\Leftrightarrow\quad (\dot{\bvarepsilon}_p,\dot{\gamma_p})\in \partial I_{\mathcal{E}}(\Sigma_E,g)=N_{\mathcal{E}}(\Sigma_E,g)\,,
\end{equation}
where $I_{\mathcal{E}}$ is the indicator function of the set $\mathcal{E}$ of admissible generalized stresses and $N_{\mathcal{E}}(\Sigma_E,g)$ is the normal cone of the set $\mathcal{E}$ at $(\Sigma_E,g)$.\\ The condition
 (\ref{equivflowlaw})$_2$ is called the {\it dual form of the flow
law}, which in the case of smoothness of the yield surface
$\partial\mathcal{E}$ at $(\Sigma_E,g)$ gives for some scalar parameter $\lambda\geq0$
\begin{equation}\label{KKT-iso}
\dot{\bvarepsilon}_p=\lambda\, \frac{\dev\Sigma_E}{\norm{\dev\Sigma_E}}\quad\mbox{ and }\quad
\dot{\gamma}_p=\lambda=\norm{\dot{\bvarepsilon}_p}\end{equation}
together with the Karush-Kuhn-Tucker complementary conditions:
$$\lambda\geq0\,, \quad\phi(\Sigma_E,g)\leq0 \quad\mbox{ and }\quad 
\lambda\,\phi(\Sigma_{\mbox{\scriptsize E}},g)=0\,.$$
Note that with this choice, the global dissipation inequality (\ref{global-reduced-dissip}) is satisfied.
\subsubsection{Mathematical strong formulation of the model}\label{strong}
Taking into account the free energy density $W$ in
(\ref{free-eng-iso}) together with the additional term in (\ref{new-term1}) and the constraint $\norm{q}\leq\beta$ in the
definition of the dissipation function $\mathcal{D}$ in
(\ref{diss-funt-iso}), the model is strongly formulated as follows:
 find
\begin{itemize}\item[(i)] the displacement $u\in \SFH^1(0,T;
\SFH^1_0(\Omega,{\Gamma},\mathbb{R}^3))$\,,
\item[(ii)] the infinitesimal plastic strain $\bvarepsilon_p\in
\SFH^1(0,T;L^2(\Omega, \sL(3)\cap\mbox{Sym}(3)))$ with
\begin{eqnarray*} \sym\Curl \bvarepsilon_p &\in &\SFH^1(0,T; L^2(\Omega,\mbox{Sym}(3)\cap\sL(3)));\\
\inc\bvarepsilon_p=\Curl[(\Curl \bvarepsilon_p)^T] &\in & \SFH^1(0,T;
L^2(\Omega,\BBR^{3\times 3}));\\
 \inc\inc\bvarepsilon_p=\Curl\Bigl(\bigl[\Curl\Curl[(\Curl \bvarepsilon_p)^T] \bigr]^T\Bigr)&\in& \SFH^1(0,T;
L^2(\Omega,\BBR^{3\times 3}))\,,\end{eqnarray*}
\item[(iii)] The internal isotropic hardening  variable $\gamma_p\in
\SFH^1(0,T;L^2(\Omega))$\,,\end{itemize} such that the content of
Table \ref{table:isohard-inc-new} holds.

\begin{table}[h!]\footnotesize\begin{center}
\begin{tabular}{|ll|}\hline &\qquad\\
 {\em Additive split of strain:} &$\sym\nabla u=\bvarepsilon_e+\bvarepsilon_p$\\
{\em Equilibrium:} & $\mbox{Div}\,\sigma +f=0$ with
$\sigma=\C_{\mbox{\scriptsize{iso}}}\bvarepsilon_e$\\&\\ {\em
Free energy:} &
$\dsize\frac12\,\langle\C_{\mbox{\scriptsize{iso}}}\bvarepsilon_e,\bvarepsilon_e\rangle+\frac12\,\mu\,
k_2\,|\gamma_p|^2 +\,\,\frac12\,\mu\,L^2_c\norm{\sym\Curl\bvarepsilon_p}^2 +\,\frac12\,\mu\,\widehat{L}^4_c\norm{\inc\bvarepsilon_p}^2$\\&\\
{\em Yield condition:} &$\phi(\Sigma_E,g):=\norm{\dev\Sigma_E}^2 -(\yieldzero-g)\leq0$\quad where \quad$g:=-\mu\,k_2\gamma_p$\\&$\begin{array}{lcl} \Sigma_E&:=&\sigma+\Sigma_{\mbox{\scriptsize
curl}} +\Sigma_{\mbox{\scriptsize
 inc}}\,,\\
  \Sigma_{\mbox{\scriptsize curl}}&:=&-\,\mu\,L^2_c\,\sym\Curl(\sym\Curl\varepsilon_p)\,\\
\Sigma_{\mbox{\scriptsize
inc}}&:=&-\mu\,\widehat{L}^4_c\,\inc\inc\bvarepsilon_p\,=\,\Curl\bigl[\Curl\Curl[\Curl \bvarepsilon_p]^T \bigr]^T\end{array}$
 \\&\\{\em Dissipation inequality:} &
 $\dsize\int_\Omega\bigl[\la\Sigma_E,\dot{\bvarepsilon}_p\ra+g\,\dot{\gamma}_p\ra\bigr]\,dx\geq0$\\
 {\em Dissipation function:} &$\mathcal{D}(\dot{\bvarepsilon}_p,\dot{\gamma}_p):=\left\{\begin{array}{ll}\yieldzero\, \norm{\dot{\bvarepsilon}_p} &\mbox{ if }
 \norm{\dot{\bvarepsilon}_p}\leq\dot{\gamma}_p,\\
 \infty &\mbox{ otherwise}\end{array}\right.$\\&\\
 {\em Flow law in primal form:} &
 $(\Sigma_E,g)\in\partial \mathcal{D}(\dot{\bvarepsilon}_p,\dot{\gamma}_p)$\\
{\em Flow law in dual form:}
&$\dot{\bvarepsilon}_p=\lambda\,\dsize\frac{\dev\Sigma_E}{\norm{\dev\Sigma_E}},
\quad\quad \dot{\gamma}_p=\lambda=\norm{\dot{\bvarepsilon}_p}$\\
{\em KKT conditions:} &$\lambda\geq0$, \quad $\phi(\Sigma_E,g)\leq0$,
\quad $\lambda\,\phi(\Sigma_E,g)=0$\\
 {\em Boundary conditions for $\bvarepsilon_p$:} & $\bvarepsilon_p\times n|_{\partial\Omega}=0\,,\quad (\Curl \bvarepsilon_p)^T\times n|_{\partial\Omega}=0$\\&\\
 \hline
\end{tabular}\caption{\footnotesize The new regularized irrotational model with isotropic hardening and Kr\"{o}ner's incompatibility tensor.  Also in this case, the boundary condition on $\bvarepsilon_p$
necessitates at least $\bvarepsilon_p,\,\,(\Curl
\bvarepsilon_p)^T\in \mbox{H}(\mbox{Curl};\,\Omega,\,\BBR^{3\times
3})$. The model is micro-random i.e., invariant w.r.t. $p\to p+A(x)$,  $A(x)\in\so(3)$ and invariant under infinitesimal conformal mappings $p\to p+\nabla\phi_c$ with $\dev\sym\nabla\phi_c=0$.}\label{table:isohard-inc-new}\end{center}\end{table}

\subsubsection{Weak formulation of the model}\label{wk-inc-iso} To
obtain the weak formulation of the model, we consider the
equilibrium in its weak formulation. That is, for every $v\in
\SFH^1_0(\Omega,\mathbb{R}^3)$ we have
\begin{equation}\label{weak-eq-iso-new}
\int_{\Omega}\la\C_{\mbox{\scriptsize{iso}}}(\sym\nabla
u-\bvarepsilon_p),\mbox{sym}(\nabla v-\nabla\dot{u})\ra
dx=\int_\Omega f(v-\dot{u})dx\,.
\end{equation}
On the other hand, for every $q\in
C^\infty(\overline{\Omega},\sL(3)\cap\mbox{Sym}(3))$ such that
$$q\times n|_{\partial\Omega}=0\quad\mbox{ and }\quad (\Curl q)^T\times n|_{\partial\Omega}=0\,$$ and for every $\beta\in
L^2(\Omega)$, integrate (\ref{primal-flow-law}) over $\Omega$ using
 the pair of functions $(q,\beta)$ and get 
\begin{eqnarray}\label{primal-weak1}
\int_\Omega\mathcal{D}(q,\beta)\,dx &\geq &\int_\Omega\mathcal{D}(\dot{\bvarepsilon}_p,\dot{\gamma}_p)\,dx +\int_\Omega\bigl[\la\sigma +\Sigma_{\mbox{\scriptsize curl}}+\Sigma_{\mbox{\scriptsize inc}},q-\dot{\bvarepsilon}_p\ra +g\,(\beta-\dot{\gamma}_p)\bigr]\,dx\\
\nonumber &=&\int_\Omega\bigl[\la\sigma -\mu\,L^2_c\,\sym\Curl(\sym\Curl\varepsilon_p)-\mu\,\widehat{L}^4_c\,\inc\inc\bvarepsilon_p,q-\dot{\bvarepsilon}_p\ra\\
&&\nonumber \hskip2truecm-\mu\,k_2\,\gamma_p\,(\beta-\dot{\gamma}_p)\bigr]\,dx\,.\end{eqnarray}
Now integrating by parts the two terms
$\la\sym\Curl(\sym\Curl\varepsilon_p),q-\dot{\bvarepsilon}_p\ra$ once and 
 $\la\inc\inc\bvarepsilon_p,q-\dot{\bvarepsilon}_p\ra$ twice, using the
boundary conditions
$$(q-\dot{\bvarepsilon}_p)\times n|_{\partial\Omega}=0\quad\mbox{ and }\quad (\Curl (q-\dot{\bvarepsilon}_p))^T\times n|_{\partial\Omega}=0\,$$ we get from (\ref{primal-weak1}) that
\begin{eqnarray}
\label{weak-dissip-ineq-new}\nonumber \int_\Omega
\mathcal{D}(q,\beta)\,dx &\geq &
\int_\Omega \mathcal{D}(\dot{\bvarepsilon}_p,\dot{\gamma}_p)\,dx +\int_\Omega\la \C_{\mbox{\scriptsize{iso}}}(\sym\nabla
u-\bvarepsilon_p),\,q-\dot{\bvarepsilon}_p\ra\,
dx\\
&&\nonumber-\mu\,k_2\int_\Omega \gamma_p(\beta-\dot{\gamma}_p)\,dx\,-\,\mu\,L^2_c\int_\Omega\la\sym\Curl
\bvarepsilon_p,\,\sym\Curl(q-\dot{\bvarepsilon}_p)\ra
\,dx
 \\&& -\,\mu\,\widehat{L}^4_c\int_\Omega\la\inc
\bvarepsilon_p,\,\inc(q-\dot{\bvarepsilon}_p)\ra
\,dx\,.
\end{eqnarray}
Adding (\ref{weak-dissip-ineq-new}) to the weak formulation of the
equilibrium in (\ref{weak-eq-iso-new}), we get the weak formulation of
 our model of gradient plasticity with isotropic hardening and
Kr\"{o}ner's incompatibility tensor

\begin{eqnarray}\label{weak-form-iso-new}
&&\nonumber\hskip-1truecm\int_{\Omega}\Bigl[\la\C_{\mbox{\scriptsize{iso}}}(\sym\nabla
u-\bvarepsilon_p),(\sym\nabla
v-q)-(\sym\nabla\dot{u}-\dot{\bvarepsilon}_p)\ra+\mu\,
k_2\gamma_p(\beta-\dot{\gamma}_p)\\
&&\nonumber\quad +\,\mu\,L^2_c\,\la\sym\Curl
\bvarepsilon_p,\,\sym\Curl(q-\dot{\bvarepsilon}_p)\ra+\,\mu\,\widehat{L}^4_c\la\inc\bvarepsilon_p,\,\inc(q-\dot{\bvarepsilon}_p)\ra\Bigr]\,dx\\
&&\quad+\,\int_\Omega \mathcal{D}(q,\beta)\,dx
-\int_\Omega \mathcal{D}(\dot{\bvarepsilon}_p,\dot{\gamma}_p)\,dx\,\geq\,\int_\Omega
f(v-\dot{u})\,dx\,.\end{eqnarray} That is,
\begin{equation}\label{VI}\ba(w,z-\dot{w})+j(z)-j(\dot{w})\geq\la\ l,z-\dot{w}\ra\,,\end{equation}
where
\begin{eqnarray}
\nonumber\ba(w,z)&:=&\int_{\Omega}\Bigl[\la\C_{\mbox{\scriptsize{iso}}}(\sym(\nabla
u-p)),\sym\nabla v-q\ra+\mu\, k_2\,\gamma_p\,\beta+ \\
&&\hskip2truecm  +\,\,\mu\,L_c^2\,\la\sym\Curl
\bvarepsilon_p,\,\sym\Curl q\ra+\,\mu\,\widehat{L}_c^4\la\inc\bvarepsilon_p,\,\inc q\ra\Bigr]dx\,, \label{bilin-form-iso}
\\\nonumber\\
 j(z)&:=&\int_\Omega
\mathcal{D}(q,\beta)\,dx\,,\label{functional-iso}\\
\langle l,z\rangle&:=&\int_\Omega
fv\,dx\,,\label{lin-form-iso}\end{eqnarray} for $w=(u,\bvarepsilon_p,\gamma_p)$ and
$z=(v,q,\beta)$\,.

\subsubsection{Existence result for the weak
formulation}\label{exist-weak-iso-new} We prove the existence result for the weak formulation
(\ref{weak-form-iso-new}) by closely following the approach by now classical, which uses the abstract machinery developed by
Han and Reddy in \cite{Han-ReddyBook} for mathematical problems in geometrically linear
classical plasticity and used for instance in \cite{DEMR1, REM, NCA,
EBONEFF, ENR2015} for models of gradient plasticity.  Precisely, we will need the following Theorem.
\begin{theorem}\label{han-reddy-theo}{\rm{\bf(\cite[Theorem 6.19]{Han-ReddyBook})}}\\
Let $\SFZ$ be a Hilbert space and let $\SFW$ be a nonempty closed convex cone in $\SFZ$. Consider the following problem: find $w\in\SFH^1([0,T]; \SFZ)$ with $w(0)=0$ such that for almost every $t\in[0,T]$, $\dot{w}(t)\in\SFW$ and
\begin{equation}\label{HR-ABS1}
\ba(w,z-\dot{w})+j(z)-j(\dot{w})\geq \langle
\ell,z-\dot{w}\rangle\mbox{ for every } z\in \SFW\,.\end{equation}
Assume that the following hold:\begin{itemize}\item[1.] the bilinear form $\ba$ is symmetric, continuous on $\SFZ$ and coercive on $\SFW$, i.e., there exist $C>0$ and $\alpha>0$ such that
\begin{equation}\label{HR-ABS2}
 \ba(w,z)\leq C\,\norm{w}_Z\,\norm{z}_Z\qquad\forall w,\,z\in\SFZ\quad\mbox{ and }\quad
\ba(z,z)\geq\alpha\,\norm{z}^2_Z\qquad\forall z\in \SFW\,;\end{equation}
\item[2.] $\ell\in \SFH^1([0,T];\SFZ')$  with $\ell(0)=0$.

\item[3.] the functional $j$ is non-negative, convex,  lower continuous  and positively $1$-homogeneous  $\SFW$, i.e., $j(sz)=|s|\,j(z)\quad\forall s\in\mathbb{R}\,,\quad\forall z\in\SFZ\,.$
\end{itemize} Then the problem (\ref{HR-ABS1}) has a solution $w\in \SFH^1([0,T];\SFZ)$.
\end{theorem}

Therefore, the problem is then reduced to finding a suitable Hilbert space $\SFZ$ and its subset $\SFW$ such that the bilinear form $\ba(\cdot,\cdot)$ and the functionals $j$ and $\ell$ satisfies the assumptions of Theorem \ref{han-reddy-theo}.
\vskip.2truecm\noindent The choices of function spaces for the displacement variable $u$ and the isotropic hardening variable $\gamma_p$ are straightforward as 
$$u\in \SFH^1_0(\Omega,\mathbb{R}^3)\quad\mbox{ and }\quad \gamma_p\in \SFL^2(\Omega)\,.$$ For the plastic strain variable $\bvarepsilon_p$, we first need to introduce the space
\begin{eqnarray}\label{space-plastic1}
\nonumber \hskip-.4truecm\SFH_{\mbox{\scriptsize inc }}(\Curl,\,\Omega;\sL(3)\cap\mbox{Sym}(3))&:=&\{q\in \SFL^2(\Omega,\,\sL(3)\cap\mbox{Sym}(3)\,\,|\,\,(\Curl q)^T\in \SFH(\Curl,\,\Omega;\mathbb{R}^{3\times 3})\}\\
 \hskip-.4truecm&:=&\{q\in \SFH(\Omega,\,\sL(3)\cap\mbox{Sym}(3)\,\,|\,\,\inc q\in \SFL^2(\Curl,\,\Omega;\mathbb{R}^{3\times 3})\}\end{eqnarray}
equipped with the norm
\begin{eqnarray}\label{norm-space-plastic2}
\nonumber \norm{q}^2_{\mbox{\scriptsize inc}}&:=&\norm{q}^2_{L^2(\Omega)}+\norm{(\Curl q)^T}^2_{\mbox{\scriptsize H}(\mbox{\scriptsize
Curl};\Omega)}\,=\,\norm{q}^2_{L^2}+\norm{(\Curl q)^T}^2_{L^2}+\norm{\Curl(\Curl q)^T}^2_{L^2}\\
&=&\norm{q}^2_{\mbox{\scriptsize H}(\mbox{\scriptsize
Curl};\Omega)} +\norm{\inc q}^2_{L^2(\Omega)} \,.\end{eqnarray}
Let us mention that spaces of functions involving the  inc-operator were already used in the literature and we refer the interested reader for instance to the papers \cite{AMSVANGOETH, AMSVANGOETH2017}.

We also consider 
the closure
$\SFH_{\mbox{\scriptsize sym, inc }}(\Curl,\,\Omega,\partial\Omega;\sL(3)\cap\mbox{Sym}(3))$ of
the linear subspace
$$\Bigl\{q\in C^\infty(\overline{\Omega},\mbox{Sym}(3))\,\,|\,\,\tr{q}=0,\,\,q\times\,n|_{\partial\Omega}=0\mbox{ and }(\Curl q)^T\times n|_{\partial\Omega}=0\Bigr\}$$ in the
norm
\begin{equation}\label{newnorm-q}\norm{q}_{\mbox{\scriptsize
symcurl,\,inc}}^2:=\norm{q}^2_{L^2}+\norm{\sym\mbox{Curl\,}q}^2_{L^2}+\norm{\inc q}^2_{L^2}\,.
\end{equation}
 Motivated by the well-posedness question for models of infinitesimal gradient plasticity (specially for models dictated by invariance under infinitesimal rotations) \cite{NESNEF2012, NESNEF2013,EBONEFF,NSW2009,NCA}, infinitesimal Cosserat elasticity \cite{NJ2009,JN2010, NEFF2006}, infinitesimal Cosserat elasto-plasticity \cite{NC2005,NK2008,CN2008,NCMW2007} and infinitesimal relaxed micromorphic \cite{NGMP2014,NJMR2009,NEFFFIN2006}, Bauer et al. \cite{BNPS2014,BNPS2015} (see also Neff et al. \cite{NPW2011-1, NPW2012-1, NPW2012-2, NPW2014}) 
 derived a new
inequality extending Korn's first inequality to incompatible tensor
fields, namely there exists a constant $C(\Omega)>0$ such that
\begin{align}
\label{incompatible_korn}
\forall \, X\in \SFH(\Curl;\,\Omega,\,\BBR^{3\times 3}), \quad & X\times\,n|_{\partial\Omega}=0:   \\
  & \|X\|_{L^2(\Omega)}\le C(\Omega)\,
     \Big( \norm{\sym X}_{L^2(\Omega)}+
      \norm{\Curl X}_{L^2(\Omega)}\Big)\, .\notag
\end{align}
Now, if we apply the incompatible Korn's type inequality to $X=(\Curl q)^T$ for $q\in  C^\infty(\overline{\Omega},\mbox{Sym}(3))$ with $(\Curl q)^T\times n|_{\partial\Omega}=0$, we get
\begin{eqnarray}\label{equi-est1}
\nonumber\norm{\Curl q}_{L^2(\Omega)} &=&\norm{(\Curl q)^T}_{L^2(\Omega)}\\
\nonumber&\leq& C(\Omega)\Bigl(\norm{\sym\Bigl((\Curl q)^T\Bigr)}_{L^2(\Omega)} +\norm{\Curl[(\Curl\bvarepsilon_p)^T]}_{L^2(\Omega)}\Bigr)\\
&=&C(\Omega)\Bigl(\norm{\sym\Curl q}_{L^2(\Omega)} +\norm{\inc q}_{L^2(\Omega)}\Bigr)\,,\end{eqnarray}
then we have the decisive identity
\begin{eqnarray}\label{space-eps}
\nonumber&&\hskip-1.5truecm\SFH_{\mbox{\scriptsize sym, inc
}}(\Curl,\,\Omega,\partial\Omega;\sL(3)\cap\mbox{Sym}(3))\\
&\hskip.5truecm\equiv&\bigl\{q\in
\SFH_0(\Curl,\,\Omega,\partial\Omega;\sL(3)\cap\mbox{Sym}(3)),\quad (\Curl q)^T\in \SFH_0(\Curl,\,\Omega,\partial\Omega;\mathbb{R}^{3\times 3})\bigr\}\\
\nonumber &\hskip.5truecm=& \bigl\{q,\,(\Curl q)^T\in
\SFH(\Curl,\,\Omega,\mathbb{R}^{3\times 3}),\quad \tr q=0\mbox{ a.e. in }\Omega,\quad q\times n|_{\partial\Omega}=(\Curl q)^T\times n|_{\partial\Omega}=0\bigr\}\end{eqnarray}with the norms $\norm{\cdot}_{\mbox{\scriptsize inc}}$ and $\norm{\cdot}_{\mbox{\scriptsize 
symcurl,\,inc}}$ being equivalent.
\vskip.2truecm\noindent
 Now,
we set
\begin{eqnarray}
\mbox{V}:&=&\mbox{H}^1_0(\Omega;\,\mathbb{R}^3)\,,\label{space-u-iso}\\
 \mbox{Q}:&=&\bigl\{q\in \SFH_{\mbox{\scriptsize inc }}(\Curl,\,\Omega;\sL(3)\cap\mbox{Sym}(3))\,\,|\,\,q\times n|_{\partial\Omega}=0\mbox{ and }(\Curl q)^T\times n|_{\partial\Omega}=0\bigr\}\,,\label{space-plastic}\\
 \Lambda:&=& \mbox{L}^2(\Omega)\,,\label{space-intvar}\\
   \mbox{Z}:&=&\mbox{V}\times \mbox{Q}\times\Lambda\,,\label{product-space-iso}\\
   \mbox{W}:&=&\{z=(v,q,\beta)\in\mbox{Z}\,\,|\,\,
    \norm{q}\leq\beta\}\,,\label{conv-set-iso}
\end{eqnarray} equipped with the norms
\begin{eqnarray}
&&\nonumber \norm{v}_V:=\norm{\nabla
v}_{L^2(\Omega)},\label{norm-V-iso}\quad\qquad
\norm{q}_{Q}:=\norm{q}_{\mbox{\scriptsize inc}}\,\,,\label{norm-Q-iso}\quad\qquad\norm{\beta}_{\Lambda}=\norm{\beta}_{L^2(\Omega)}\,,\label{norm-Lambda-iso}\\
&&\norm{z}^2_{Z}:=\norm{v}^2_{V} +\norm{q}^2_{Q}
+\norm{\beta}^2_{\Lambda}\quad\mbox{ for } z=(v,q,\beta)\in Z\,.
\label{norm-Z-iso}
\end{eqnarray}
Let us prove the coercivity of the bilinear form $\ba(\cdot,\cdot)$ on the closed
convex set $\SFW$, where the constraint $\norm{q}\leq\beta$ in $\SFW$ plays
a crucial role. Let therefore $z=(v,q,\beta)\in \SFW$. Then,
\begin{eqnarray*}
\ba(z,z) &\geq&  m_0\norm{\sym(\nabla v)-
q}^2_2+\mu\,k_2\norm{\beta}^2_{L^2} +\mu\, L^2_c\,\norm{\sym\Curl q}^2+\mu\,
\widehat{L}_c^4\norm{\inc q}_{L^2}^2\\
&&\hskip3truecm \mbox{ ($m_0>0$ is from (\ref{ellipticityC}))}\\
&=&m_0\left[\norm{\sym (\nabla v)}^2_2 +\norm{q}^2_{L^2} -2\la\sym
(\nabla v), p\ra\right]+\mu\,k\norm{\beta}_{L^2}^2 \\
&&\hskip3truecm +\mu\, L^2_c\,\norm{\sym\Curl q}^2+\mu\, \widehat{L}_c^4\,\norm{\inc q)}_{L^2}^2
\\
&\geq &m_0\left[\norm{\sym(\nabla v)}^2_{L^2} +\norm{ q}^2_{L^2}
-\theta\norm{\sym (\nabla
v)}_{L^2}^2-\frac1\theta\norm{q}_{L^2}^2\right] +\frac12\mu\,k_2\norm{\beta}_2^2\\
&& +\frac12\mu\,k_2\norm{ q}_{L^2}^2 + \mu\,L^2_c\,\norm{\sym\Curl
q}^2_{L^2} +\mu\,\widehat{L}_c^4\norm{\inc q}_{L^2}^2
\\
&&\mbox{ (using Young's inequality and
 }\norm{q}\leq\beta\mbox{ from }W)\\
&=& m_0(1-\theta)\norm{\sym (\nabla
v)}^2_2+\left[m_0(1-\frac1\theta)+\frac12\mu\,k_2\right]\norm{
q}_2^2 +\frac12\mu\,k_2\norm{\beta}_2^2\\&& +\mu\,L_c^2\norm{\sym\Curl q}^2_{L^2}+
\mu\,\widehat{L}_c^4\norm{\inc q}_{L^2}^2\,.\end{eqnarray*} So, choosing $\theta$ such that
$\displaystyle\frac{2\,m_0}{2\,m_0+ \mu\,k_2}\leq \theta<1,$ and
using the classical Korn's first inequality, there exists some positive constant
$K(m_0,\,\mu\,,k_2\,,\Omega)>0$ such that
\begin{eqnarray}\label{ineq-coerc1}\ba(z,z)&\geq& K\left[\norm{v}_V^2+\norm{q}^2_{L^2}
+\norm{\beta}_\Lambda^2 + \mu\,L^2_c\,\norm{\sym\Curl
q}^2_{L^2} +\mu\,\widehat{L}_c^4\norm{\inc q}_{L^2}^2 \right]\\
\nonumber &\geq & C\left[\norm{v}_V^2+\norm{q}^2_{Q}
+\norm{\beta}_\Lambda^2\right]\,=\,C\norm{z}^2_{Z}\quad\forall
z=(v,q,\beta)\in W\,,\end{eqnarray}
where $C=C(m_0,\,\mu\,,k_2\,,\Omega, L_c,\widehat{L}_c)>0$. 
For the second inequality in (\ref{ineq-coerc1}), we used the inequality (\ref{equi-est1}) obtained as a consequence of Korn's type inequality for incompatible tensor fields in Neff et al. 
\cite{NPW2011-1, NPW2012-1, NPW2012-2, NPW2014}\,.\hfill\qed

\vskip.2truecm\noindent So, assuming that the body is initially unloaded and undeformed, which corresponds to assuming that $f(x,0)=0$ for almost all $x\in\Omega$ with homogeneous initial conditions, we obtained the following existence result for the weak formulation (\ref{weak-form-iso-new}) of our model.
\begin{theorem}\label{existence-wf-iso-new} Under the choices of the Hilbert space $\SFZ$ and the closed convex cone $\SFW$ in (\ref{space-u-iso})-(\ref{conv-set-iso}) with the norms in (\ref{norm-Z-iso}) and the functionals $\ba$, $j$ and
$\ell$ in (\ref{bilin-form-iso})-(\ref{lin-form-iso}), the weak formulation (\ref{weak-form-iso-new}) when written as the variational inequality of the second kind (\ref{HR-ABS1}) has a solution $w=
(u,\bvarepsilon_p,\gamma_p)$ in $ \SFH^1([0,T];\SFZ)$ with $\dot{w}\in\SFL^2([0,T];\SFW)$.\end{theorem}

\begin{remark}\label{uniqueness-rem}{\rm Uniqueness of the strong solution is obtained as in \cite{EHN2016} provided the following further assumptions are satisfied:
\begin{eqnarray}\label{reg-assump}
\sym\Curl\bigl(\sym\Curl\bvarepsilon_p\bigr) &\in& \SFL^2(\Omega,\mbox{Sym}(3)\cap\sL(3))\,\\
\nonumber \Curl\Bigl([\Curl\Curl[\Curl\bvarepsilon_p]^T]^T\Bigr) &\in&\SFL^2(\Omega,\mathbb{R}^{3\times 3})\,.\end{eqnarray}

}\end{remark}
\section{Discussion}
It remains a difficult task to reconcile mathematical and physical requirements. Indeed, the incorporation of  Kr\"{o}ner's incompatibility tensor $\inc\bvarepsilon_p$ is physically transparent and the novel model is micro-random and gauge-invariant.  Micro-randomness being useful for polycrystals and gauge-invariance being a generally physically necessary requirement. However, using integration by parts in order to arrive at a global reduced dissipation inequality, the following  lowest order boundary conditions 
\begin{equation}\label{bc-lowest1}
\bvarepsilon_p\times n|_{\partial\Omega}=0\quad\mbox{ and }\quad (\Curl\bvarepsilon_p)^T\times n|_{\partial\Omega}=0\end{equation} impose themselves. 

From a mathematical point of view these expressions are, however, not well-defined as boundary traces through  a  control of the given free-energy. In order to give them a well-defined meaning, we resorted to adding an additional term in the free-energy, namely $$\frac12\,\mu\,L^2_c\norm{\dev\sym(\Curl\bvarepsilon_p)^T}^2=\frac12\,\mu\,L^2_c\norm{\sym(\Curl\bvarepsilon_p)^T}^2\,,.$$ the equality here is due to the fact that $\tr(\Curl S)=0$ for ever  This term provides the missing boundary control for (\ref{bc-lowest1}) by  Korn's-type inequality for incompatible tensor fields in Neff et al. \cite{NPW2011-1, NPW2012-1, NPW2012-2, NPW2014,BNPS2014, BNPS2015}\,. However, the additional term breaks the gauge-invariance of the model, while it satisfies the micro-randomness condition.

On the positive side, the invariance under the diffeomorphism group (gauge-invariance) is replaced by the invariance under infinitesimal  conformal group (both statements adapted to our geometrically linear setting).

At the moment, we do not know how to set up a theory which is fully gauge-invariant and micro-random, while at the same time being mathematically well-posed. Consider e.g. a model with plastic spin and add $\frac12\,\mu\,L_c^2\,\norm{\Curl p}^2$ (see Table \ref{table:isohard-fullyiso-new}). This choice does not provide any control of $\norm{\Curl\sym p}^2=\norm{\Curl\bvarepsilon_p}^2$ necessary for well-posedness of (\ref{bc-lowest1}).

A preliminary conclusion could be that the micro-randomness assumption, which effectively reduces the flow law to the six-dimensional space of symmetric plastic strains $\bvarepsilon_p$, is to be critically seen in gradient plasticity approaches which are also supposed to  satisfy gauge-invariance.

\addcontentsline{toc}{section}{Acknowledgements}

\section*{Acknowledgements:} 
We thank David J. Steigmann (University of Berkeley) for inspiring discussions on invariance conditions in plasticity theory which motivated us to introduce our micro-randomness condition and  to consider the fourth order gradient plasticity model for polycrystals. The first author thanks the Faculty of Mathematics of the University of Duisburg-Essen (Germany) for its kind hospitality during his visit  in April 2017.

\addcontentsline{toc}{section}{Appendix}
\section*{Appendix}\label{appendix}
{ Let us first establish that 
\begin{equation}\label{trace-curl-sym1} \tr\Curl S=0\qquad\forall S:\mathbb{R}^3\to\mbox{Sym}(3) \mbox{ smooth tensor field}.\end{equation}

In fact, recalling that $(\Curl S)_{ij}=\sum_{kl=1}^3\epsilon_{ikl}\,S_{jl,k}$ we have
\begin{eqnarray}\label{trace-curl-sym2}\nonumber\tr(\Curl S)&=&\sum_{i}(\Curl S)_{ii}=\sum_{i}\sum_{kl}\epsilon_{ikl}\,S_{jl,k}\\
&=&(S_{13,2}-S_{12.3})+(S_{21,3}-S_{23,1})+(S_{32,1}-S_{31,2})=0\,,
\end{eqnarray} because  $S\in\mbox{Sym}(3)$ and hence, $S_{13,2}=S_{31,2}$, $S_{12,3}=S_{21,3}$ and $S_{23,1}=S_{32,1}$.
\vskip.3truecm
Below are some further 
properties of the Kr\"{o}ner's incompatibility tensor defined by
\begin{equation}\label{Kroner}
\inc(\bvarepsilon_p):=\Curl[(\Curl\bvarepsilon_p)^T]\,.
\end{equation} For the convenience of the reader,
 we note that
\begin{eqnarray}&&\inc(\bvarepsilon_p)\in\mbox{Sym}(3)\quad\mbox{ since
}\quad\bvarepsilon_p\in\mbox{Sym}(3)\,;\label{inc-id1}\\
&&\inc(\inc(\bvarepsilon_p))\in\mbox{Sym}(3)\,;\label{inc-id2}\\
&&\mbox{tr}(\inc(\bvarepsilon_p))=\Delta\,\mbox{tr}(\bvarepsilon_p)
-\div(\Div\bvarepsilon_p)=-\div(\Div\bvarepsilon_p)\quad\mbox{since
}\bvarepsilon_p\in\sL(3)\,;\label{inc-id3}\\
&&\mbox{tr}(\inc(\inc
\bvarepsilon_p)=-\Delta\div(\Div\bvarepsilon_p)\,.\label{inc-id4}\end{eqnarray}
Since (\ref{inc-id2}) follows from (\ref{inc-id1}), let us establish
here the identities (\ref{inc-id1})-(\ref{inc-id4}) for the reader's
convenience. First of all, in components
$$(\Curl\bvarepsilon_p)_{ij}:=\epsilon_{jkl}\,(\bvarepsilon_p)_{il,k}\quad \Leftrightarrow\quad (\Curl\bvarepsilon_p)^T_{ij}:=
\epsilon_{ikl}\,(\bvarepsilon_p)_{jl,k}\,.$$
 Hence
\begin{equation}\label{inc-component1}(\inc{\bvarepsilon_p})_{ij}=\Curl[(\Curl\bvarepsilon_p)^T]_{ij}=\epsilon_{jkl}\,(\Curl\bvarepsilon_p)^T_{il,k}
=\epsilon_{jkl}\,\epsilon_{imn}\,(\bvarepsilon_p)_{ln,mk}\,.\end{equation}
Now, notice that $(\bvarepsilon_p)_{ln,mk}=(\bvarepsilon_p)_{nl,km}$.
Therefore,
\begin{equation}\label{inc-component2}(\inc{\bvarepsilon_p})_{ij}=\epsilon_{jkl}\,\epsilon_{imn}\,(\bvarepsilon_p)_{ln,mk}
=\epsilon_{jmn}\,\epsilon_{ikl}\,(\bvarepsilon_p)_{nl,km}=
\epsilon_{jmn}\,\epsilon_{ikl}\,(\bvarepsilon_p)_{ln,mk}=(\inc(\bvarepsilon_p))_{ji}\,,\end{equation}
which establishes (\ref{inc-id1}). Now,
\begin{equation}\label{trace-inc}\mbox{tr}(\inc(\bvarepsilon_p))=(\inc(\bvarepsilon_p))_{ii}=\epsilon_{ikl}\,\epsilon_{imn}\,
(\bvarepsilon_p)_{ln,mk}\,.\end{equation} Using the identity
$$\epsilon_{ikl}\,\epsilon_{imn}=\delta_{km}\,\delta_{ln}-\delta_{kn}\,\delta_{lm}\,
\quad\mbox{ (see for instance, \cite[epsilon-delta identity
(1.20)]{GURTAN-BOOK}})\,,$$ we get
\begin{eqnarray}\label{inc-component3}
\nonumber\mbox{tr}(\inc(\bvarepsilon_p))&=&\delta_{km}\,\delta_{ln}\,(\bvarepsilon_p)_{ln,mk}-\delta_{kn}\,\delta_{lm}\,(\bvarepsilon_p)_{ln,mk}
\,=\,(\bvarepsilon_p)_{ll,mm}
-(\bvarepsilon_p)_{lk,lk}\\
\nonumber&=&(\bvarepsilon_p)_{ll,mm}
-(\bvarepsilon_p)_{kl,lk}\,=\,\Delta\mbox{tr}(\bvarepsilon_p)-(\Div\bvarepsilon_p)_{k,k}\\
&=&\Delta\mbox{tr}(\bvarepsilon_p)-\div(\Div\bvarepsilon_p)=-\div(\Div\bvarepsilon_p)\,,\end{eqnarray}
which establishes (\ref{inc-id3}). So, from (\ref{inc-id3}), it
follows that if $\bvarepsilon_p$ is a divergence-free tensor or
$\Div\bvarepsilon_p$ is a divergence-free vector field, then
$\inc(\bvarepsilon_p)$ becomes trace-free, that is,
$\inc(\bvarepsilon_p)\in\sL(3)$. Now, to establish (\ref{inc-id4}),
notice that
$$\Div\Curl X=0 \quad\forall
X\in C^2(\Omega,\,\BBR^{3\times 3})\,.$$ This trivially follows from
our definitions of $\Curl$ and $\Div$ of a second tensor field as
row-wise operations. Hence, $\Div(\inc\bvarepsilon_p)=0$. So, using
(\ref{inc-id3}), we find that
\begin{eqnarray}
\nonumber
\mbox{tr}(\inc(\inc\bvarepsilon_p))&=&\Delta\mbox{tr}(\inc\bvarepsilon_p)
-\div(\Div(\inc\bvarepsilon_p))=\Delta^2(\mbox{tr}(\bvarepsilon_p)-\Delta(\div(\Div\bvarepsilon_p))\\
&=&-\Delta(\div(\Div\bvarepsilon_p))=-\div(\Div\Delta\bvarepsilon_p)\,,\end{eqnarray}
where $\Delta^2=\Delta(\Delta)$ denotes the bi-Laplacian operator.
Therefore, the tensor $\inc(\inc\bvarepsilon_p)$ is trace-free if
one of the conditions below satisfied:
\begin{itemize}\item[(i)] $\bvarepsilon_p$ is a divergence-free
tensor field;
\item[(ii)]$\Div\bvarepsilon_p$ is a divergence-free vector field;\item[(iii)]
$\div(\Div\bvarepsilon_p)$ is an harmonic function;\item[(iv)]
$\bvarepsilon_p$ is an harmonic tensor field;
\item[(v)] $\Delta\bvarepsilon_p$ is a divergence-free tensor field;
\item[(vi)] $\Div(\Delta\bvarepsilon_p)$ is a divergence-free vector
field.\end{itemize}

}

\end{document}